\documentclass[smallextended]{svjour3}       
\usepackage{graphicx}        
\usepackage{algorithm2e}
\usepackage{amsmath}
\usepackage{amsfonts}
\usepackage{multirow}
\usepackage{tabularx}
\usepackage{array}
\usepackage{float} 
\usepackage{placeins}
\begin{document}

\title{Locally refined spline surfaces for terrain and sea bed data: tailored approximation,
export and analysis tools}
\titlerunning{Representing terrain and sea bed with LR B-spline surfaces}
\author{Vibeke Skytt \and Tor Dokken}

\institute{V. Skytt \\ Corresponding author \at
  SINTEF Digital, PO Box 124 Blindern, 0314 Oslo, Norway \\
  \email{Vibeke.Skytt@sintef.no} \\
  Orcid: 0000-0001-5117-958X
  \and
  T. Dokken
   \at
  SINTEF Digital, PO Box 124 Blindern, 0314 Oslo, Norway \\
  \email{Tor.Dokken@sintef.no} \\
  Orcid: 0000-0001-6000-9961
}

\maketitle

\begin{abstract}
The novel Locally Refined B-spline (LR B-spline) surface format is   suited for representing terrain and seabed data in a compact way. It provides an alternative to the well know raster and triangulated surface representations. An LR B-spline surface has an overall smooth behaviour and allows the modelling of local details with only a limited growth in  data volume. In regions where many data points belong to the same smooth area LR B-splines allow a very lean representation of the shape by locally adapting the resolution of the spline space to the size and local shape variations of the region.
The surfaces generated  approximate the smooth component of a cloud of data points within user specified tolerances. The method can be  modified to improve the accuracy in particular domains and selected data points. The resulting surfaces are well suited for analysis and computing secondary information such as contour curves and minimum and maximum points. The surfaces can be translated into a raster  or a tessellated surface of desired quality or exported as collections of tensor product spline surfaces. Data transfer can also be performed  using Part 42 of  ISO 10303 (the STEP standard)  where LR B-splines were published in 2018. 
\end{abstract}

{\noindent
{\bf Funding}
This work has been supported by the Norwegian research counsil under grant number 270922.
}

{\noindent
{\bf Acknowledgements}
  The Norwegian map authorities, division Sj\o kartverket
has provided the measurement data.
}

\section{Introduction} \label{sec:intro}
Geospatial  data acquisition  of terrains and sea beds produces huge point
clouds. The structure, or lack of structure, in the point cloud depends on the actual data acquisition technology used. Efficient downstream uses of the acquired data require in general  structured and compact data representations. Working directly in the point cloud is not attractive.

 The approach we  follow in this paper is to use the novel  Locally Refined B-spline surfaces (LR B-splines)
 for representation of geospatial data. To relate LR B-splines to surface representations already in use in Geographic Information Systems (GIS),
  we provide in Section~\ref{sec:GISsurf} an overview of other relevant representations.
 
 The starting point for an LR B-spline surface is  a tensor product B-spline surface.  In the LR B-spline surface  extra degrees of freedom are inserted locally when needed. The desired  smoothness is maintained, and  the growth in data volume is kept under control. The extension of tensor product B-spline surfaces  
to LR B-splines is presented in detail in Section~\ref{sec:LR}. A summary and comparison of the presented surface representation formats is given in Section~\ref{sec:compareformats}. Section~\ref{sec:LRApprox} presents the procedure for creating LR B-spline surfaces from geographical point clouds. 
Section~\ref{sec:LRApprox} also include an example where we
compare approximation with LR B-spline surfaces and raster
representations of the same point cloud.

Interoperability with existing GIS-technology is essential for the deployment of LR B-splines within GIS. Section~\ref{export} addresses how LR B-spline surfaces can be converted to standard formats used in GIS-systems. 

Once an LR B-spline surface representation of some GIS data is obtained, this surface can
be interrogated and analyzed to compute and derive terrain  information, e.g., contour curves, slope and aspect ratio. In Section~\ref{analysis} we focus on contour curves
and extremal points.

The format of locally refined splines is flexible. The property of local
refinement is used to expand the modelling freedom in areas with high
degree of local variation in the data set. This flexibility can be exploited
further by requiring higher accuracy in specified areas  or in selected points. In the subsea context it is most often important to have high accuracy in areas with
shallow water rather than where the water is  deep. Section~\ref{adopt_apprx} looks into three 
different methods to emphasize particular domains in the data set. 

Finally, 
we will provide some conclusions and set directions for further work in
Section~\ref{conclusion}.

 \section{GIS surface representations} \label{sec:GISsurf} 
 Representation of terrain phenomena has been a topic 
 for a long time and  not less so in the digital era. Different aspects with digital terrain modelling
 is discussed in~\cite{terrain}. Among the topics are representation of digital terrain model surfaces and quality
 assessment.
 
 Spatial interpolation is the process of using points with known values to estimate values at unknown sample points. 
 The corresponding interpolator can be exact or approximative depending on whether or not the initial data points are fitted exactly. In the context of large geospatial data sets only approximative interpolators are feasible.
 Methods using all data points to estimate a sample point are global, while methods using a subset of points are local.
 A survey on interpolation methods can be found in~\cite{GISinterpolation}. 
 
Most GIS systems offer the possibility of the Inverse distance weighting interpolation (IDW) method
originally proposed by Shepard, see~\cite{IDW}.
The method is simple, but tends to create flat 
spots around sparse data points. Initially IDW was a global method, but modifications have been
proposed to define sample points from subsets of
data and to avoid flat spots~\cite{AltShepard,ShepardLim}.

Natural neighbour interpolation~\cite{NNneighbour} is originally based
on Voronoi tessellation and tends to give a more smooth
result than the Inverse distance weighting.
Kriging~\cite{Kriging} is an advanced technique using a Geo-statistical approach where the estimation errors are directly available. Kriging is a weighted average
technique where all data points are analyzed to find the auto-correlation of the
particular data set being interpolated. Radial basis functions (RBF) provide another set of methods for scattered data interpolation. Originally a global method resulting in equation systems of the same size as the number of data points but local methods have also been developed~\cite{LocalRBF}.  Several types of basis functions are used in the context of RBF, and Franke evaluated of some functions in~\cite{RBFTest}. Buhmann discusses
radial basis functions and properties thereof in~\cite{radial}.

In contrast to the previously presented methods
the use of splines in GIS interpolation can result in sample values outside the convex hull of the data points. This restrains the smoothing effect inherent in those methods and enables the estimation of lows and
highs not present in the input data, but may result
in exaggerated behaviour, in particular if the data points are close and have extreme differences in value. Splines in GIS most often mean splines in tension or regularized splines, see~\cite{GISinterpolation}, which differ from the tensor product splines presented in Section~\ref{sec:Bcurve}.

 The raster representation is the most frequently used data format in GIS. To make a raster representation of a  point cloud, spatial interpolation is used to define values, for instance elevation or precipitation, in a regular grid. 
 The digital elevation model (DEM), which is a raster representation, is the
 most commonly used format for processed terrain and sea bed data.
 The raster is a compact, highly structured and efficient representation. It is also an approximate representation as the input data are not exactly fitted. The accuracy depends on the resolution of the raster and the selected interpolation method. If the resolution is low compared to the variation of the data then the result will be inaccurate, if it is high 
 then the data volume grows more than necessary. If there are large  differences in the local variation of the data in different areas, then a trade-off between accuracy and data volume must be made. 
 
 A raster represented surface is not completely defined by the
 sample points. Given a raster the interpolation method is often unknown, and values between the sample values have to be estimated. Several approaches are available. The value in the cell center may be selected. Another simple method is 
 bivariate evaluation where the estimated value is computed from a
 bivariate surface interpolating the four surrounding sample values.
 Inverse distance weighting and other interpolation methods can be used also in this context.

Root mean square error (RMSE) is the standard way of measuring the 
error of a model. $RMSE = \sqrt\frac{\sum_{i=1}^n(\hat{z}_i-z_i)^2}{n}$,
where $n$ is the number of data values, $z_i$ are the
given values and $\hat{z}_i$ are the corresponding estimated
values. RMSE is a good measure if we want to estimate the standard deviation $\sigma$ of a typical observed value from our model’s prediction, assuming that our observed data can be decomposed into the predicted value and random noise with mean zero.
It does not, however, take the spatial pattern of the error into account. Fisher and Tate
discuss errors in digital elevation models in~\cite{ErrDEM1}. Different interpolation
methods are considered without giving one clear recommendation. The accuracy of which
a given terrain is approximated with a DEM depends on the match between the DEM
resolution and the terrain features. Wise investigates the relation between terrain
types and DEM error in~\cite{ErrDEM2} and evaluates five interpolation methods (bilinear, IDW, RBF, spline and local polynomial) in the context of re-sampling in~\cite{ErrDEM3}. The
effect on derivative information such as slope is also studied. Splines and RBF is
found to give the best results. He also finds a link between the characteristics of the
terrain and the pattern of errors in elevation.

A triangulated irregular network (TIN) is a continuous surface representation  frequently used in GIS. This is a flexible format for geospatial data that allows adaptation to local variations. It is thus capable of very high accuracy. In~\cite{TINbasin} a triangulated surface is used to represent a drainage-basins while hydrological similarity is used in the TIN creation in~\cite{TINhydro}. To restrict the data size, approximation is required also here. However, the nodes of a TIN are distributed variably to create an accurate representation of the terrain. TINs can, thus, have a higher resolution in areas where a surface is highly variable or where more detail is desired and a lower resolution in areas that are less variable.

TIN models have a more complex data structure than raster surface models and tend to be more expensive to build and process. 
Therefore TINs are typically used for high-precision modelling of smaller areas. The triangulated surface
consists of flat triangles. Points in-between the corner points in a triangle are found by linear interpolation, which can give a jagged appearance of the surface.
The problem is especially visible at sharp or nearly sharp edges, but can be 
remedied by methods like constrained Delaunay triangulation. 
However, conforming Delaunay triangulation is often recommended over constrained triangulation
due to less propensity to create long, skinny triangles, which are undesirable for surface analysis. In~\cite{terrain} a comprehensive discussion on various aspects with triangulation in terrain modelling is presented.

\section{B-splines and locally refined splines} \label{sec:LR}
Within Computer-Aided Design (CAD) the univariate minimal support B-splines basis is used for the representation of B-spline spline curves. For the representation of tensor product B-spline surfaces in CAD a basis generated by the tensor product of two univariate minimal support B-splines bases is used.
In the context of CAD the term "B-spline basis function" is often used, this is valid for the minimal support bases described above. However, when building locally refined spline spaces for T-splines and Locally refined B-splines we risk, if care is not taken, to end up with a collection of tensor product B-splines spanning the spline space, that is not linearly independent. Consequently the collection of tensor product B-splines that span the spline space will not always constitute a basis for the B-spline space spanned. To explain Locally Refined B-splines we need to recap a number of B-spline basics. In Section \ref{sec:bsplines} we address a single B-spline, the tensor product of two B-splines, the univariate minimal support B-spline bases and the bivariate minimal support tensor product  B-spline basis. Then we address B-spline curves and tensor product  B-splines surfaces in  Section~\ref{sec:Bcurve}, and locally refined B-spline surfaces in~\ref{sec:LRBsurf}.

\subsection{B-splines, tensor product B-splines and minimal support B-spline bases}\label{sec:bsplines}
Given a non-decreasing sequence $\mathbf u=(u_0,u_1,\ldots,u_{p+1})$ we define a B-spline
$B[\mathbf u]:{\mathbb{R}} \to {\mathbb{R}}$ of degree $p\ge 0$ recursively as follows
\begin{equation}
B[\mathbf u](u):=\frac{u-u_0}{u_{p}-u_0} B[u_0,\ldots,u_{p}](u)+
\frac{u_{p+1}-u}{u_{p+1}-u_1}B[u_1,\ldots,u_{p+1}](u),\label{eq:Bsp-local-knots}
\end{equation}
starting with
$$B[u_i,u_{i+1}](u):=\begin{cases} 1;&\text{if $u_i\le u < u_{i+1}$};\\
0;&\text{otherwise},\end{cases}
\quad i=0,\ldots,p.
$$
We define $B[\mathbf u]\equiv 0$ if $u_{p+1}=u_0$ and
terms with zero denominator are defined to be zero.

A univariate spline space can be defined by a polynomial degree $p$ and a knot vector $\mathbf u = \{u_0, u_2, \ldots, u_{N+p} \}$, where the knots  satisfy: $u_{i+1} \ge u_{i}$, $i=0,\ldots ,N+p-1$, and  $u_{i+p+1} > u_{i}$, $i=0,\ldots,N-1$. I.e, a knot value can be repeated $p+1$ times. The number of times a knot value is repeated is called the multiplicity $m$ of the knot value. The continuity across a knot value of multiplicity $m$ is $C^{p-m}$.

A basis for the univariate spline space above can be defined in many ways. However, the approach most often used is the  univariate minimal support B-spline basis. In this the B-splines are  defined by selecting $p+2$ consecutive knots from $\mathbf u$, starting from the first knot. So $B_{i,p}(u):= B[u_i,\ldots,u_{i+p+1}](u)$ is defined by the knots $u_i,\ldots,u_{i+p+1}$, $i=0,\ldots,N-1$ . The minimal support B-spline basis has useful properties such
as local support, non-negativity and partition of unity
(the basis functions sum up to one in all parameter values in the interval $[u_p,u_{N}]$).

Given two non-decreasing knot sequences $\mathbf u = \{u_0, u_1, \ldots, u_{N_1+p_1} \}$ and $\mathbf v = \{v_0, v_1, \ldots, v_{N_2+p_2} \}$ where  $p_1\ge0$ and $p_2\ge0$. We define a bivariate
tensor product B-spline $B_{i,j,p_1,p_2}: {\mathbb{R}^2} \to {\mathbb{R}}$ from the two univariate B-splines 
$B_{i,p_1}(u)$ and $B_{j,p_2}(v)$ by
$$B_{i,j,p_1,p_2}(u,v):=B_{i,p_1}(u) B_{j,p_2}(v).$$
The support of $B$ is given by the cartesian product
$$\text{\em supp}(B_{i,j,p_1,p_2}):=[u_i,u_{i+p_1+1}]
\times[v_j,v_{j+p_2+1}].$$

A  bivariate tensor product spline space is made by  the tensor product of two univariate spline spaces. Assuming that both univariate spline spaces have a minimal support B-spline basis, the minimal support basis for the tensor product B-spline space is constructed by making all tensor product combinations of the B-splines of the two bases. The minimal support basis for this  spline space  contains the tensor product B-splines $B_{i,j,p_1,p_2}(u,v)$, $i=0,\ldots,N_1-1$, $j=0,\ldots,N_2-1$. As in the univariate case the basis  has useful properties such as 
non-negativity and partition of unity.

\subsection{B-spline curves and tensor product B-spline surfaces} \label{sec:Bcurve}
Spline curves are frequently represented using a univariate minimal support B-spline basis. 
\begin{equation}
f(u) = \sum_{i=0}^{N-1} c_i B_{i,p}(u), u \in [u_{p},u_{N}].
\end{equation}
Here  $c_i \in R^d$, $i=0,\ldots,N-1$ are the curve coefficients and $d$ is the dimension of the geometry space.  
The curve lies in the convex hull of its coefficients.

\begin{figure}
\centering
\includegraphics[width=0.9\textwidth]{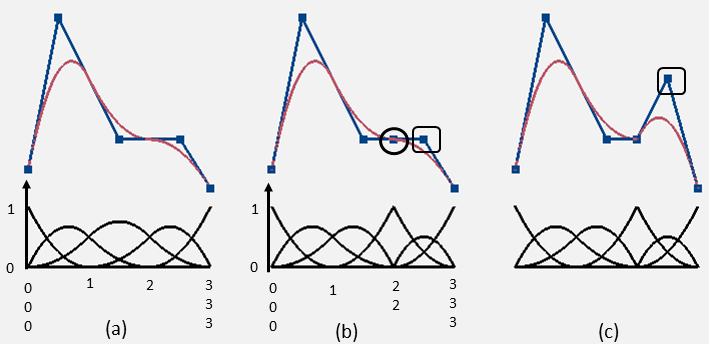} 
\caption{Knot insertion into a quadratic B-spline curve. (a) Initial curve with basis functions,
  (b) curve and basis functions after knot insertion, (c) curve and basis functions after re-positioning of one coefficient}
\label{fig:Bsplinecurve}
\end{figure}
A B-spline curve can be locally refined. Fig.~\ref{fig:Bsplinecurve} (a) shows a quadratic curve with knot vector ${\{0,0,0,1,2,3,3,3\}}$. The curve coefficients and the control polygon corresponding to the curve are included in the 
figure, and the associated B-splines
are shown below. In (b)  a new knot with value $2$ is added, thus
increasing the knot multiplicity in an already existing knot. The curve is not altered, but
the control polygon is enhanced with a new coefficient interpolating the curve. The coefficient is marked with a circle
in Fig.~\ref{fig:Bsplinecurve} (b). The double knot at $2$ allows creating a curve with $C^0$ continuity. When moving the coefficient marked with
a square, we obtain a sharp corner (c). Note that only the last part of the curve is
modified.

Spline surfaces are frequently represented using a bivariate minimal support tensor product B-spline basis.
\begin{equation}
F(u,v)=
\sum_{i=0}^{N_1-1} \sum_{j=0}^{N_2-1} c_{i,j} B_{i,j,p_1,p_2}(u,v)  \in   [u_{p_1},u_{N_1}]\times[v_{p_2}, v_{N_2}].
\label{eq:TPsurface}
\end{equation}
Here  $c_{i,j} \in R^d$, $i=0,\ldots,N_1-1$, $j=0,\ldots,N_2-1$ are the surface coefficients and $d$ is the dimension of the geometry space. 

\begin{figure}
\centering
\includegraphics[width=0.6\textwidth]{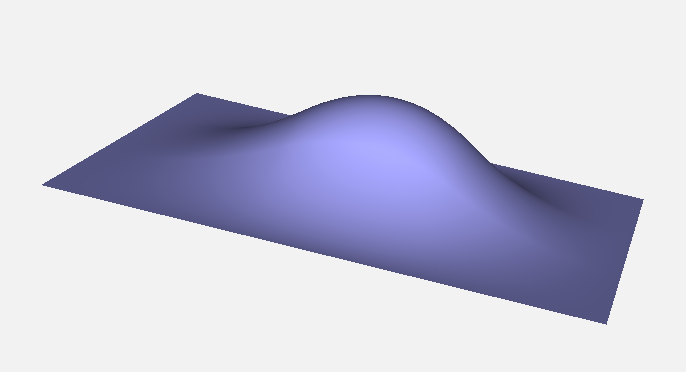} \\
\caption{Example of a bi-quadratic tensor product B-spline basis functions.}
\label{fig:B-splines}
\end{figure}

The construction of the basis functions implies that the properties
of the univariate B-spline basis functions carry over to the tensor product case, implying that the tensor product spline surface
is a stable construction. Fig.~\ref{fig:B-splines} shows an example of a tensor product B-spline with the knot vector 
${\{0,1,2,3\}}$ in both parameter directions.

Tensor product B-spline surfaces lack the property of local refinement.
The distance between consecutive knots may differ so the size of the polynomial
pieces may vary and increased knot multiplicity decreases the continuity between
adjacent polynomial segments. However, a knotline cannot be restricted to only
a part of the surface. If one new knot is inserted in the first parameter direction of the surface  (\ref{eq:TPsurface})  the number of coefficients will increase with $N_2$.

\subsection{Locally refined B-spline surfaces} \label{sec:LRBsurf}
In some applications such as representation of terrain and sea bed, the lack of
local refinement is a severe restriction. A tensor product surface covers
a rectangular domain and the need for approximation power will not be uniformly
distributed throughout this domain. There are three main approaches for extending
spline surfaces to support local refinement compared 
in~\cite{trivariate}.
\begin{itemize}
    \item {\bf Hierarchical B-splines}~\cite{hierarchical} are based on a dyadic sequence of grids determined by scaled lattices. On each level of the dynastic grids a spline space spanned by uniform tensor product B-spines is defined. The refinement is performed one level at the time. Tensor product B-splines on the coarser level are removed and B-splines at the finer level added in such a way that linear independence is guaranteed. The sequence of spline space for Hierarchical B-splines will be nested.
    \item {\bf T-splines}~\cite{tsplines}  denote a class of locally refined splines, most often presented using bi-degree $(3,3)$. The starting point for T-spline refinement is a  tensor product B-spline surface with control points and meshlines (initial T-mesh) with assigned knot values. For a bi-degree (3,3) T-splines the knot vectors of the tensor product B-spline corresponding to a control point is determined by moving  in the T-mesh outwards from the control point in all four axis parallel parameter directions and picking in each direction knot values from the two first T-mesh lines intersected. The mid knot value in each knot vector is copied from the control point. The refinement is performed  by successively adding new control points in-between two adjacent control points in the T-mesh. The control point inherits one parameter value from the T-mesh line. The other parameter is chosen so that sequence of control point parameter values have a monotone evolution along the T-mesh line. Control points on adjacent parallel T-mesh lines that have one shared parameter value are connected with a new T-mesh line. 
    The general formulation of T-splines does not guarantee a sequence of nested spline spaces and that the polynomial space is spanned over  each element. However, T-spline subtypes such as semi-standard T-splines and Analysis Suitable 
    T-splines \cite{AST} do.
    \item {\bf Locally Refined B-splines}~\cite{lrsplines}, the refinement approach of this paper, starts (as T-splines) from a tensor product B-spline surface. The refinement is performed successively by inserting axis parallel meshlines in the mesh of knotlines. Each meshline inserted has to split the support of at at least one tensor product B-spline. The constant knot value of a meshline inserted is used for performing univariate knot insertion. The refinement is performed in the parameter direction orthogonal to the meshline in all tensor product B-splines that have a support split by the meshline. This approach ensures that the spline spaces produced are nested and that the polynomial space is spanned over all polynomial elements. In Section \ref{sec:LRref} we provide further details on additional refinements triggered and how the resulting tensor product B-splines can be scaled to achieve partition of unity. 
\end{itemize}

\subsubsection{Preparing the tensor product B-spline surface for local refinement} \label{sec:TPlocalref}
After the first local refinement is performed on the initial tensor product B-splines surface the regular mesh structure is lost. The representation in  Equation \ref{eq:TPsurface}, where the local knot vectors are described by an index pointing into a global knot array and the polynomial degree, is no longer applicable. We now organize the  tensor product B-splines, each now knowing its local knot vectors, as a collection of B-splines $\mathcal{B}_0$. Doing this we can represent the tensor product surface  in Equation \ref{eq:TPsurface}  by setting $\mathcal{B}_j=\mathcal{B}_0$.
\begin{equation}
F(u,v)=\sum_{B \in \mathcal{B}_j} c_B s_B B(u,v).\label{eq:tensor-collection}
\end{equation}
Here $c_B$ is the coefficient corresponding to the B-spline $B$.
We have added a scaling factor $s_B$ to each tensor product B-spline. This is used for recording  how the  collection of tensor product B-splines generated during the refinement have to be scaled to form a positive partition of unity. For a tensor product B-spline surface all $s_B \equiv 1$. Using the representation we can remove B-splines to be refined by knot insertion from the collection and replace them with the refined B-splines and a scaling factor that provides a scaled partition of unity. Note that the same tensor product B-spline can result from the refinement of different B-splines. Consequently duplicates have to be identified and their scaling factors accumulated.

\subsubsection{LR B-spline refinement}\label{sec:LRref}

The process of Locally Refine B-splines (LR B-splines) is described in detail in \cite{lrsplines}.   Please consult the paper for formal proofs. Below a summary of the most important steps is presented. The refinement always starts from a tensor product B-spline space represented as in Equation \ref{eq:tensor-collection}. The refinement proceeds  with a sequence of meshline insertions producing a series of collections of tensor product B-splines $\mathcal{B}_0, \mathcal{B}_1, \ldots, \mathcal{B}_j, \mathcal{B}_{j+1}$ spanning nested spline spaces. Note  that we require that all tensor product B-splines in these collection have minimal support. By this we mean that all meshlines that cross the support of a tensor product B-spline also have to be a knotline of the tensor product B-splines counting multiplicity. Thus, to ensure that all tensor product B-splines have this property the process of going from $\mathcal{B}_j$ to $\mathcal{B}_{j+1}$ frequently includes a number of steps and a sequence of LR B-spline collections $\mathcal{B}_{j,1},\mathcal{B}_{j,2},\ldots,\mathcal{B}_{j,r_j}$, where $\mathcal{B}_{j,1}=\mathcal{B}_j$ and $\mathcal{B}_{j,r_j}=\mathcal{B}_{j+1}$.

Assume that Equation \ref{eq:tensor-collection} represents $F(u,v)$ using the collection $\mathcal{B}_{j}$ and we now want to represent $F(u,v)$ after refining with a given  meshline to produce a refined representation of $F(u,v)$. The refinement process goes as follows:

\begin{itemize}
    \item As long as there is a $B\in \mathcal{B}_j$ that does not have minimal support on the refined mesh we proceed as follows.
Let $\gamma$ be a meshline that splits the support of $B$. This means that either $\gamma$ is not a knotline of $B$, or $\gamma$ is a knotline of $B$ but has higher multiplicity than the knotline of $B$. Decompose $B$ into its univariate component B-splines $B(u,v)=B(u)B(v)$. We now have two cases:
\begin{itemize}
    \item If  $\gamma$ is parallel to  second parameter direction then it has a constant parameter value $a$ in the first parameter direction. We insert $a$ in  the univariate B-spline $B(u)$ using Equation \ref{eq:boehmt} below and express $B(u)$ as $B(u)=\alpha_1 B_1(u)+\alpha_2 B_2(u)$. Then we make  two new tensor product B-splines $B_1(u,v)= B_1(u)B(v)$ and $B_2(u,v)= B_2(u)B(v)$. 
   \item If  $\gamma$ is parallel to  first parameter direction then it has a constant parameter value $a$ in the second parameter direction. We insert $a$ in  the univariate B-spline $B(v)$ using Equation \ref{eq:boehmt} below and express $B(v)$ as $B(v)=\alpha_1 B_1(v)+\alpha_2 B_2(v)$. Then we make the two new tensor product B-splines $B_1(u,v)= B(u)B_1(v)$ and $B_2(u,v)= B(u)B_2(v)$. 
\end{itemize}
$B$ can now be decomposed as follows, $B(u,v) = \alpha_1 B_1(u,v) + \alpha_2 B_2(u,v)$. We can now express $F(u,v)$ by replacing $B(u,v)$ with the two new tensor product B-splines. $F(u,v)=F(u,v) - c_B s_B B(u,v) + c_B s_B(\alpha_1 B_1(x,y) + \alpha_2 B_2(x,y))$. We update $\mathcal{B}_{j}$ by removing $B$ and adding the  tensor product B-splines $B_1(u,v)$ and $B_2(u,v)$. In addition we have to create/update both coefficients and scaling factors belonging to these two  tensor product B-splines. We must have in mind that $B_1(x,y)$ and $B_2(x,y)$ often will be duplicates of B-splines already in $\mathcal{B}_{j}$. Now let $B_r, r=1,2$
\begin{itemize}
    \item In the case $B_r$ has no duplicate set $s_{B_r}=s_B \alpha_r$ and $c_{B_r}=c_{B}$. 
    \item In the case $B_r$ has a duplicate $B_d$ set $s_{B_r} = s_{B_d}+s_B \alpha_r$ and  $c_{B_r}=({s_{B_d} c_{B_d}+s_B \alpha_r c_{B}})/ ({s_{B_r}})$, and remove the duplicate. 
   \end{itemize} 
   Note that $s_{B}$, $\alpha_r$ and $s_{B_d}$  are all positive numbers, thus $s_{B_r}$ will be positive.
   
When all $B\in \mathcal{B}_{j}$ have minimal support  we set $\mathcal{B}_{j+1}=\mathcal{B}_{j}$.
\end{itemize}
We now define the scaled tensor product B-splines $N_{B}(u,v)=s_B B(u,v)$ to provide a basis that is partition of unity for the spline space spanned by $\mathcal{B}_j$. If $F(x,y)\equiv 1$ then all coefficients of the tensor product B-spline surface we start from are $1$. In this case the coefficients $c_{B_r}$ calculated above all remain $1$, duplicates or not. Consequently,
\begin{equation}
\sum_{B \in \mathcal{B}_j}  N_B(u,v) =\sum_{B \in \mathcal{B}_j}  s_B B(u,v) =\sum_{B \in \mathcal{B}_0}  B(u,v) =1.
\label{eq:LRSurf-partition}
\end{equation}
.

 The process of inserting a knot $a \in (u_0,u_{p+1})$ into the local knot vector  $[{\mathbf{u}}]=\{u_0,\ldots,u_{p+1}\}$ belonging to a univariate B-spline $B[{\mathbf{u}}]$,  of degree $p$ was first described by Boehm~\cite{Boehm}. We organize the resulting sequence of knots as a non-decreasing knot sequence  $\{\hat{u}_0,\ldots,\hat{u}_{p+2}\}$. From this we make two new B-splines $B_1[{\mathbf u}_1]$ and $B_2 [{\mathbf u}_2] $ with corresponding knot vectors $[{\mathbf u}_1]=\{\hat{u}_0,\ldots,\hat{u}_{p+1}\}$ and $[{\mathbf u}_2]=\{\hat{u}_1,\ldots,\hat{u}_{p+2}\}$. Then
\begin{equation}
B[{\mathbf{u}}]=\alpha_1 B_1[{\mathbf u}_1] +\alpha_2 B_2[{\mathbf u}_2], \label{eq:boehmt}
\end{equation}
where 
\begin{equation}
\begin{aligned}
\alpha_1&:=\begin{cases}
 \frac{a-u_{0}}{u_{p}-u_{0}},& u_{0}<a<u_{p},\\
 1,&u_{p}\le a<u_{p+1},
\end{cases}\\
\alpha_2&:=\begin{cases}
1,&u_{0}<a\le u_{1},\\ \frac{u_{p+1}-a}{u_{p+1}-u_{1}},& u_{1}<a<u_{p+1}.
\end{cases}
\end{aligned}  
\end{equation}

\begin{figure}
\centering
\includegraphics[width=0.5\textwidth]{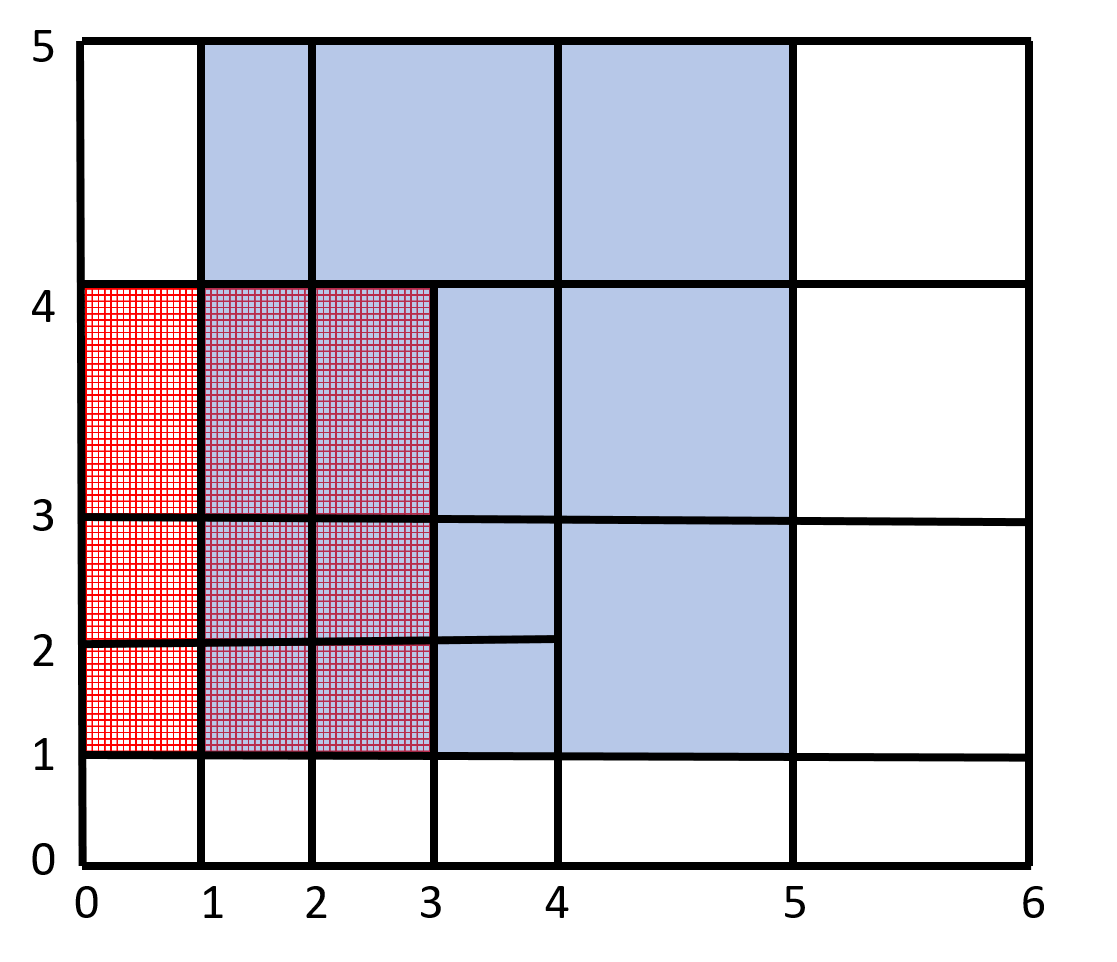} 
\caption{Parameter domain of an LR B-spline surface with indication on
  B-spline support. The knot line segments are shown as black lines. The support of two overlapping B-splines are shown in red and in blue}
\label{fig:domain}
\end{figure}
The incremental refinement by knot insertion used by Locally Refined B-splines ensures that the splines spaces generated are nested.  
Figure~\ref{fig:domain} shows a parameter
domain and the segmentation into boxes corresponding to a bi-quadratic LR B-spline surface. In addition the support of two tensor product B-splines is shown. We see that 
a knot line segment stops inside the blue tensor product B-spline support. This knot line is excluded from the local knot vectors defining the tensor product B-spline covering this support.

\subsubsection{Notation and summary} \label{sec:LRBsurf2}
To simplify notation will from now on use  the scaled tensor product B-splines defined in Equation \ref{eq:LRSurf-partition}. An LR B-spline surface of bi-degree $(p_1,p_2) $ is given as
\begin{equation}
F(u,v) =\sum_{B \in \mathcal{B}_j} c_B N_B(u,v).
\label{eq:LRSurf}
\end{equation}
The scaled tensor product B-splines used for describing an LR B-spline surface are non-negative by construction and has a compact support with size depending on the knot configuration. Partition of unity
is obtained as described in Equation \ref{eq:LRSurf-partition}. In many situations $s_B \equiv 1$, but occasionally when more tensor product B-splines than $(p_1+1)(p_2+1)$ overlap an element, then $0<s_B<1$, for some of the B-splines.

It is practical to store all knot values used in  global knot vectors to allow indexing of  knot values as in done in the representation of locally refined splines in Part 42 of ISO 10030 (STEP).
$\mathbf u_G = \{u_0, u_1, \ldots, u_{N_1+p_1} \}$ and 
$\mathbf v_G = \{v_0, v_1, \ldots, v_{N_2+p_2} \}$ as in the tensor product
case. 

Linear independence of the resulting collection $\mathcal{B}_j$ of LR B-splines is not guaranteed in general, but is
ensured for the constructions used in this article. Violations to the
linear independence can be detected and resolved by dedicated knot line
insertions. Linear dependence of LR B-splines is addressed in detail in \cite{Patrizi-2020}.

\section{Surface formats for terrain and sea bed} \label{sec:compareformats}

\begin{table}
\caption{Summary of surface formats for representing terrain and sea bed.}
{\small
\renewcommand{\arraystretch}{2}
\label{tab:DEMrepresent}
\begin{tabularx}{\textwidth} { >{\raggedright\arraybackslash}X   >{\raggedright\arraybackslash}X   >{\raggedright\arraybackslash}X   >{\raggedright\arraybackslash}X  >{\raggedright\arraybackslash}X }  
\hline
  & \bf{Representation and data structure} &\bf{Algorithm and control of accuracy} & \bf{Surface smoothness} 
& \bf{Restricting data volume} \\
\hline
\bf{Raster} & Values on regular mesh.& Spatial interpolation to define sample values. Accuracy checked after creation. & Depends on interpolation method for  evaluation between mesh points.
& Pre-set mesh resolution defining data volume.\\
\bf{Tensor product B-spline surface} & Piecewise polynomials on regular mesh, any bi-degree.
& Coefficients calculated by local/global approximation. Accuracy checked after surface creation. & Depends on polynomial degree. & Pre-set mesh resolution defining data volume. \\
\bf{TIN} & Triangulation. & Triangulate point cloud + thinning or adaptive triangulation. Accuracy can be checked during creation. & $C^0$ &  Pre-set approximation tolerance 
and/or max allowed data volume. \\
\bf{LR B-spline surface} & Piecewise polynomial on locally refined axes
parallel mesh, any bi-degree. & Coefficients calculate by  local/global approximation. Local adaption by checking accuracy during construction and refinement where needed. & Depend on polynomial degrees. & 
Pre-set approximation tolerance or restrictions in adaptive algorithm.  \\
\hline
\end{tabularx}
}
\end{table}
Table~\ref{tab:DEMrepresent} provides an overview over the surface representations
discussed in the Sections~\ref{sec:GISsurf} and~\ref{sec:LR}. TIN and LR B-spline surfaces are created by adaptive algorithms, the algorithm for LR B-spline surfaces
are presented in Section~\ref{sec:LRApprox}. Thus, degrees of freedom in the surface can be defined according to needs, and the accuracy of the
fit can be made directly available. The raster representation and the tensor product (TP) B-spline surfaces are less flexible although adaptive 
refinement can be applied to TP-surfaces as well. The lack of local refinement,
however, imply that the TP-surface size grows much faster than in the LR B-spline
case. The spline methods provide smoother surfaces than the alternatives.

Tensor product B-spline surfaces can be regarded as a generalization of the raster representation. A regular grid of polynomial elements describes the surface instead of the regular grid of points in a raster representation. In contrast
to the raster representation  the
tensor product surface is uniquely defined by its representation. Point and derivative evaluation in the 
entire surface domain is performed without ambiguity. Furthermore,
the tensor product surface description is more flexible than the raster due to the
choice of polynomial degrees, variable distances between consecutive knots and 
the possibility of defining the surface in a geometry space with dimension
higher than one.
Due to the low flexibility  the raster representation is more compact than using a tensor product B-spline surface for the representation of a piecewise constant function.

LR B-spline surfaces provides local refinement of tensor product B-spline
surfaces and the raster representation used in GIS. The flexibility is still lower than the case for TIN, but it is a smooth representation that avoids the ragged
appearance that can occur for TIN. The LR B-spline format can represent the smooth component of 
terrain and sea bed more effectively than a triangulation. Whether LR B-spline
surfaces or TIN is the most appropriate representation will depend on the
characteristics of the data and the purpose of the surface generation.
\FloatBarrier
\section{Scattered data approximation with LR B-spline surfaces}\label{sec:LRApprox}
\begin{algorithm}
\KwData {Point cloud, maximum number of iterations, threshold}
\KwResult{Approximating surface, information on approximation accuracy}
Generate initial surface\;
Compute accuracy\;
\While{there exists points with larger distance than the given threshold
and the maximum number of iterations is not reached}
{
  Refine the surface in areas where the tolerance threshold is not reached\;
  Perform approximation in the current spline space\;
  Compute accuracy\;
}
\caption{Iterative algorithm for LR B-spline surface generation}
\label{alg:approx}
\end{algorithm}

Approximation of a point cloud with an LR B-spline surface is performed in
the iterative process presented in Algorithm~\ref{alg:approx} and described in some detail
in~\cite{lralg1} and~\cite{lralg2}. We will here provide a summary. The accuracy of a current approximating surface is a measure on the
distance between the point cloud and the surface in the height direction. The maximum distance, the
average distance and the number of points, where the distance is larger than the
given threshold, is computed at every iteration step and used to
decide where the surface needs to be refined. 

In this paper we let
the x- and y-coordinates of the points serve as the parameter values while
the surface (or function) approximates the z-component of the data, but the algorithm
handles parameterized points as well (each point is given with a parameter
pair and 3D coordinates).   We will, in the
following, focus on elevation and therefore change the notation of
the surface parameters from $u$ and $v$ to $x$ and $y$.

Initially,
the point cloud is approximated by a tensor product spline surface. A tensor product B-spline surface is a special case of an LR B-spline surface and this initial surface is represented as an LR B-spline surface to start the iteration. The distances
between the surface and the points are computed in every point, and the derived accuracy information
is distributed to the polynomial elements of the surface. 
As long as the maximum distance between the points in an element and the
surface is not below the given threshold,
the element is a candidate for splitting. A broad discussion
on strategies for refinement can be found in~\cite{refineLR}. A new knot line segment must completely traverse the support of at least one B-spline as explained in Section~\ref{sec:LRref}. One new segment will increase the
number of coefficients in the surface with at least one.
To take advantage of the new degrees of freedom the surface is updated either by a least squares approximation or
by applying the multi-level B-spline algorithm (MBA)~\cite{mba}. These steps are
repeated until either the maximum distance between the surface and the point cloud
is less than the given threshold or the maximum number of iteration steps is reached.

A geospatial point cloud may represent a very
rough terrain and contain noise and possible outliers. Thus, it is not necessarily
beneficial to continue the iterations until all points are closer to the surface than the threshold. Typically, the process is stopped by the number of iterations.
In~\cite{refineLR} it is shown that the accuracy improvement tends to drastically slow down when most points have a distance to the surface less than the threshold.
Recognized outlier points may be removed during the iterative approximation process,
but distinguishing between features and outliers is not a simple task. The
topic is not pursued in this paper.

Least squares approximation is a global method where the following expression
is minimized with respect to the surface coefficients 
${c_B}$ over the entire surface domain:
$$
\min_{\mathbf c} [\alpha _1 J(F(x,y) + \alpha _2 \sum _{k=1}^K (F(x_k,y_k) - z_k)^2].
$$
$F(x,y)$ is the current LR B-spline surface, see also Equation~\ref{eq:LRSurf}, and
${\mathbf x}=(x_k,y_k,z_k), k=1,\ldots,K$
are the data points. The expression is differentiated and turned into a linear,
sparse equation system in the number of surface coefficients. 
A pure least squares approximation method will result in a singular equation
system if there exists B-splines with no data points in its support. As an
LR B-spline surface is defined on a rectangular domain and a typical point cloud
has a non-rectangular outline and may contain holes, parts of the
surface domain will frequently  lie outside the domain of the point cloud.
$J(F(x,y)$ is a smoothness term that enables a non-singular equation
system even if this situation occurs. In our case the smoothness term is an
approximation to the minimization of curvature and variation of curvature in the
surface. These originally intrinsic measures are made parameter dependent to give rise to a 
linear equation system after differentiation. The smoothness term is expressed as
$$
J(F(x,y)) = \int \int _{\Omega} \int_0^{\pi} \sum_{i=2}^3 w_i \Big(\frac{\partial ^i F(x_0+r \cos \phi,y_0 + r \sin \phi)}{\partial r^i}\Big|_{r=0} \Big) \mathrm{d} \phi \mathrm{d} x_0 \mathrm{d} y_0.
$$

At each point $(x_0,y_0)$ in the surface domain, $\Omega$, a weighted sum of the directional first and second derivatives 
of the surface is integrated around the circle and the result is integrated
over the surface domain. The directional derivative is represented
by the polar coordinates $\phi$ and $r$.The two terms in $J(F(x,y)$ are given equal weight.
A more detailed description can be found in~\cite{smooth1}. The weight on the smoothness term is kept low to emphasize
the approximation accuracy. In the examples of the next sections 
$\alpha _1= 1.0\times 10^{-9}$ and $\alpha _2 = 1-\alpha _1$. A suite of other smoothness terms is presented in~\cite{smooth2}.

Multi-level B-spline approximation (MBA) is a local, iterative approximation method, see~\cite{mba}.  
The coefficients $q_T$ of the residual surface are computed individually for each scaled tensor product B-spline $N_T$. They are determined  from the collection of data points $P=(x_P,y_P,z_P) \in \mathcal{P}_{T}$ in the support of $N_T$ as well as  the magnitude of the residuals, $r_P = z_P - F(x_P,y_P), P \in \mathcal{P}_{T}$.

\begin{equation}
q_T = \frac{\sum_{P \in \mathcal{P}_T} N_T(x_P,y_P)^2 \phi_{T,P}}{\sum _{P \in \mathcal{P}_T} N_T(x_P,y_P)^2},
\label{eq:MBA}
\end{equation}
where $\phi_{T,P}$, addressed below, depends both on the residuals of the points in $\mathcal{P}_{T}$ and  the collection of B-splines $\mathcal{B}_T$ with a support overlapping at least one point in $\mathcal{P}_{T}$. For each of the residuals we define an
under-determined equation system
$$r_P=\sum_{B\in \mathcal{B_T}} \phi_{B,P}N_{B}(x_P,y_P), P\in \mathcal{P}_T.$$
where $\phi_{B,P}$ are unknowns to be determined. There are many solutions of this equation system. For MBA  the choice is
$$
\phi_{B,P} = \frac{N_{B}(x_P,y_P)r_P}{\sum_{B'\in \mathcal{B_T}} N_{B'}(x_P,y_P)^2}, \, {B\in \mathcal{B}_T, P\in \mathcal{P}_T}.
$$
Now we  select $\phi_{T,P}, P \in \mathcal{P}_T$ to add the missing piece in Equation~\ref{eq:MBA}.
The process is explained in more detail in~\cite{mba2}.

The updated surface is the sum of the last surface and the obtained residual surface. MBA is an iterative process and repeated applications of
the algorithm converge towards the best approximation
of the point cloud given a collection of B-splines. In the  examples following MBA will be run
twice for each step in the approximation algorithm.

In~\cite{lralg1} the two approximation algorithms are compared
in a number of examples. In general, the least squares algorithm has a
better approximation order while the MBA-algorithm is more stable when the
spline space is less regular and/or the number of points in each element is
low. Typically, we use least squares approximation in the first few 
iteration steps, then turn to MBA. 

\subsection{Example: Approximating a geospatial point cloud by an LR B-spline surface} \label{sec:example}

\begin{figure}
\centering
\begin{tabular}{ccc}
\includegraphics[width=0.3\textwidth]{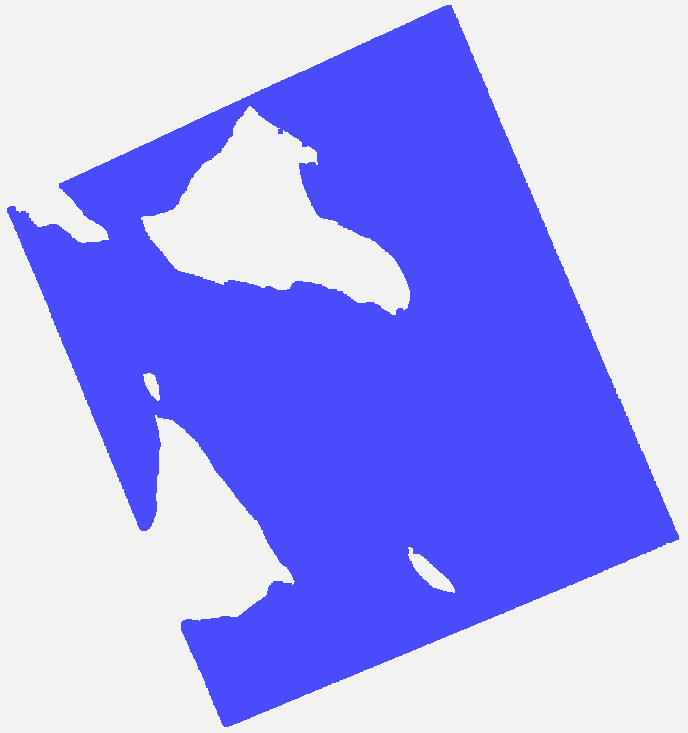}
&\includegraphics[width=0.31\textwidth]{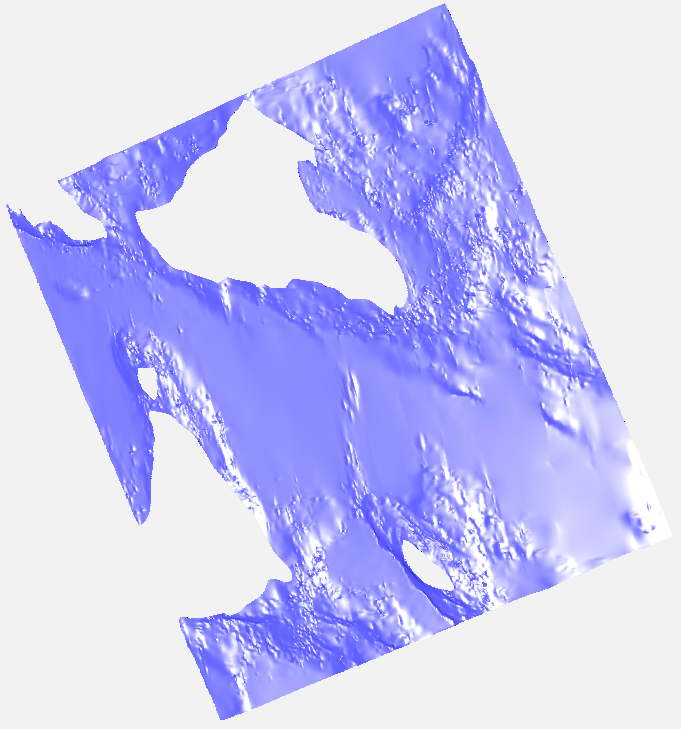}
&\includegraphics[width=0.31\textwidth]{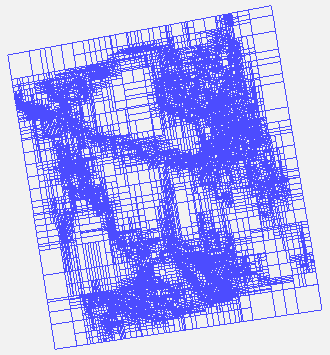}\\
(a) &(b) & (c)
\end{tabular}
\caption{(a) Outline of point cloud, (b) LR B-spline surface, (c) the element structure of the surface}
\label{fig:exs_approx}
\end{figure}

We will illustrate the LR B-spline approximation framework by applying the
algorithm to a selected data set. The same data set and corresponding surface
will be used to explain how the locally refined spline approach can be used in
a GIS work flow, and how the approximation algorithm can be tweaked to further
utilize the flexibility of the LR B-spline format.
Consider Fig.~\ref{fig:exs_approx}. The data set Fig.~\ref{fig:exs_approx}(a) covers an area of 
slightly 
less than one square kilometer close to the Norwegian coast. It consists 
of about 11 millions points with a quite uniform density although some holes.
The area is shallow, the depth varies from {-27.94} meters to {-0.55} meters. 
The point cloud contains outliers. The theoretical least maximum distance between the surface and the point cloud is 1.19 meters as the cloud contains one pair of points with the same $(x,y)$-coordinates and different depths. A surface created with 7 iterations in Algorithm~\ref{alg:approx} and
a threshold of 0.5 meters is shown in Fig.~\ref{fig:exs_approx}(b).
The data set is not compliant with a rectangular area and contains holes
while an LR B-spline surface is defined on a rectangular domain. The valid
part of the surface corresponding to this data set is defined by the use of 
trimming. The surface is bi-quadratic. This choice results in a smooth surface
with the flexibility to represent data sets with local  depth variations. The choice of
polynomial degree is investigated in~\cite{refineLR}. The sea floor covered by the data set contains both smooth areas and areas with significant variation. This can be recognized by
the element structure of the surface shown in Fig.~\ref{fig:exs_approx}(c). We can also see that
the surface is less refined in areas where there are no points.

\begin{figure}
\centering
\begin{tabular}{cc}
\multirow[b]{4}{*}[2.2cm]{\includegraphics[width=0.53\textwidth]{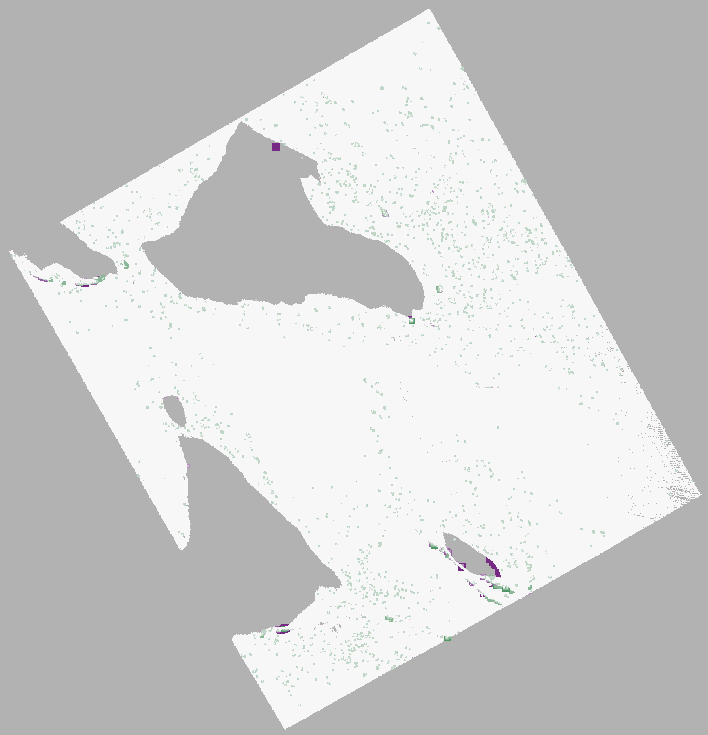}}
&\includegraphics[width=0.42\textwidth]{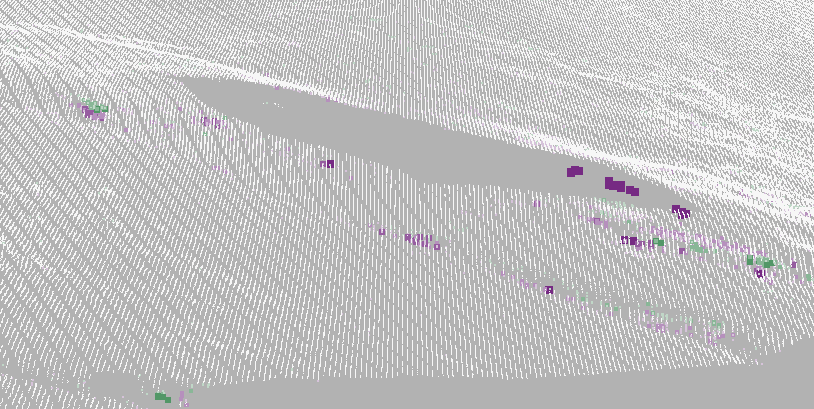} \\
&\vspace{0.1cm} \\
&\includegraphics[width=0.3\textwidth]{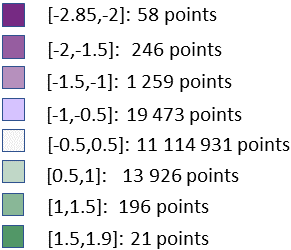} \\
&\vspace{0.05cm} \\
(a) &(b)
\end{tabular}
\caption{(a) The point cloud coloured according to the distance to the 
LR B-spline surface.
 (b) A detail and colour scale.}
\label{fig:exs_accuracy}
\end{figure}
The surface has 33\,830 coefficients and fits 99.68\% of the points within the threshold. The average distance between the point cloud and the
surface is 0.068 meters and the maximum distance is 2.85 meters. 
In Fig.~\ref{fig:exs_accuracy}(a) the point cloud is coloured 
according to the distance to the surface. Fig~\ref{fig:exs_accuracy}(b) shows a detail of the same point cloud along with the colour scale and the distance 
distribution of the points. Points that are distant to the surface
are magnified compared to closer points. The point with the maximum distance to the surface belongs to the group of points
in Fig.~\ref{fig:exs_accuracy}(c) that is slightly separated from the rest of the point cloud.

\begin{figure}
\centering
\begin{tabular}{cc}
\multirow[b]{4}{*}[2.2cm]{\includegraphics[width=0.53\textwidth]{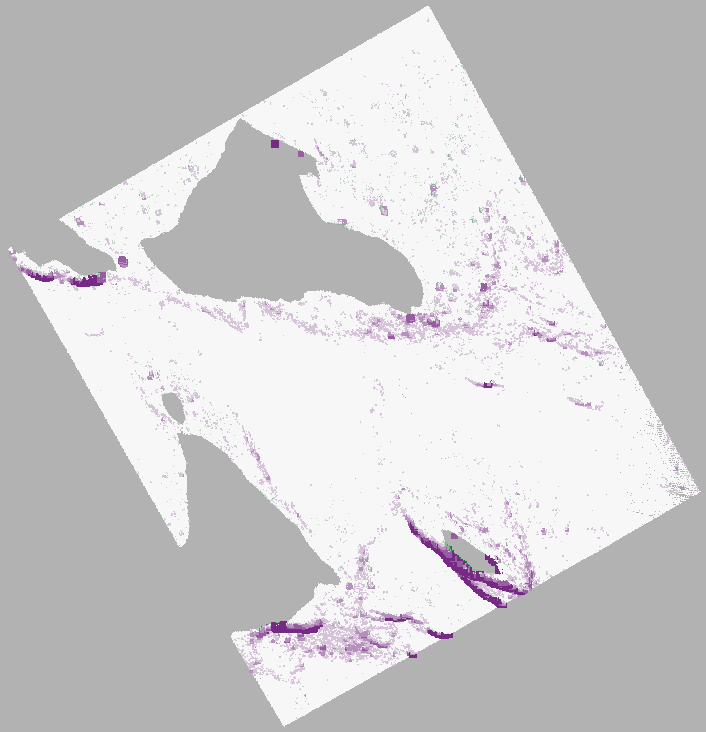}}
&\includegraphics[width=0.42\textwidth]{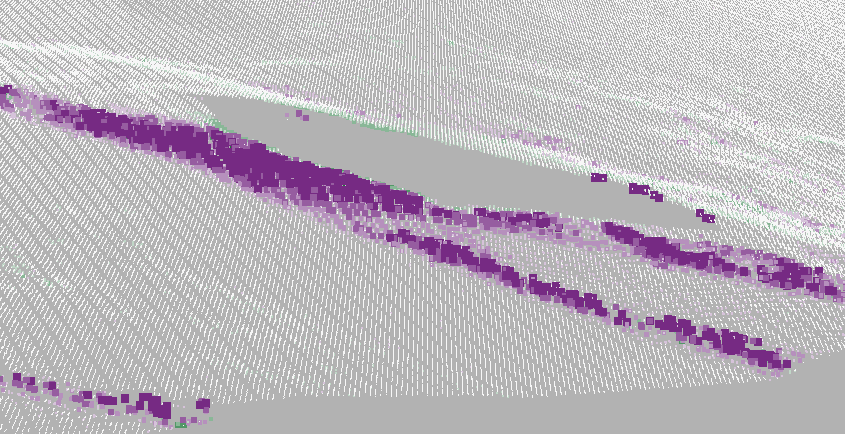} \\
&\vspace{0.1cm} \\
&\includegraphics[width=0.3\textwidth]{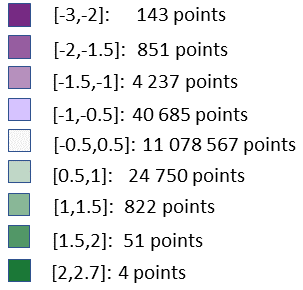} \\
&\vspace{0.05cm} \\
(a) &(b)
\end{tabular}
\caption{(a) The point cloud coloured according to the distance to the 
raster with resolution 2 meters.
 (b) A detail.}
\label{fig:raster2_accuracy}
\end{figure}
The LR B-spline surface in Fig.~\ref{fig:exs_approx} is more refined in rough than in smooth areas, but there are still more points with a
large distance to the surface in the rough areas. This effect is expected to be more dominant in a
raster representation. Fig.~\ref{fig:raster2_accuracy}(a) shows 
a set of points estimated from a $443 \times 492$ raster representation of the initial
data set. The raster has resolution 2 meters in both parameter directions and is computed with IDW. The  points displayed have the same $(x,y)$-parameters as the initial data points and the depths are estimated by
evaluating a bi-linear patch interpolating the nearby sample points. 
The points are coloured according to the distances to the initial
points and the colour coding and distance distribution is shown
in Fig.~\ref{fig:raster2_accuracy}(b) along with a detail covering the same area as in
Fig.~\ref{fig:exs_accuracy}(b). The maximum distance is 4 meters and
the average distance is 0.12 meters. The points with low accuracy
are focused in steep and ragged areas of the sea bed.

The applied inverse distance weighing interpolation uses a limited but high number of
points to estimate the sampling points with a radius, $R$, specifying the area in which points
used in the computation is collected. We set $R$ equal to 20 meters. The value $z$ of the sampling point $\mathbf{x}=(x,y)$ is
given by
$$
z(\mathbf{x}) = \frac{\sum_{i=1}^N w_i(\mathbf{x})z_i}{\sum_{i=1}^N z_i},$$
where $z_i$, $i=1,\ldots,N$ are the elevation values of the data points
used to estimate the current sampling point. The weight is 
computed as
$$
w_k(\mathbf{x}) = \Big(\frac{\max(0,R-|\mathbf{x}-\mathbf{x_k}|)}{R|\mathbf{x}-\mathbf{x_k}|}\Big)^2.$$
Here $\mathbf{x_k}=(x_k,y_k)$ is the
position corresponding to elevation value $z_k$.

\begin{figure}
\centering
\begin{tabular}{cc}
\multirow[b]{4}{*}[2.2cm]{\includegraphics[width=0.53\textwidth]{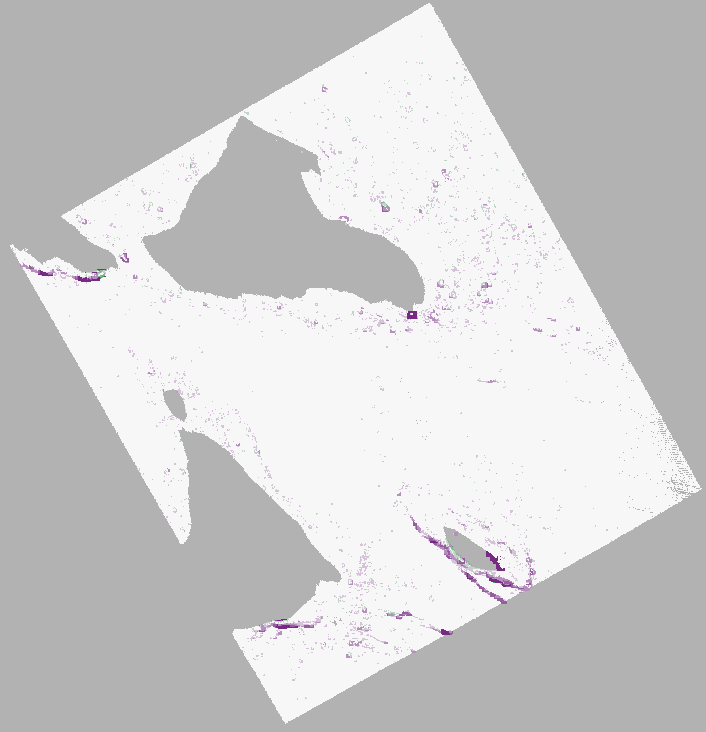}}
&\includegraphics[width=0.42\textwidth]{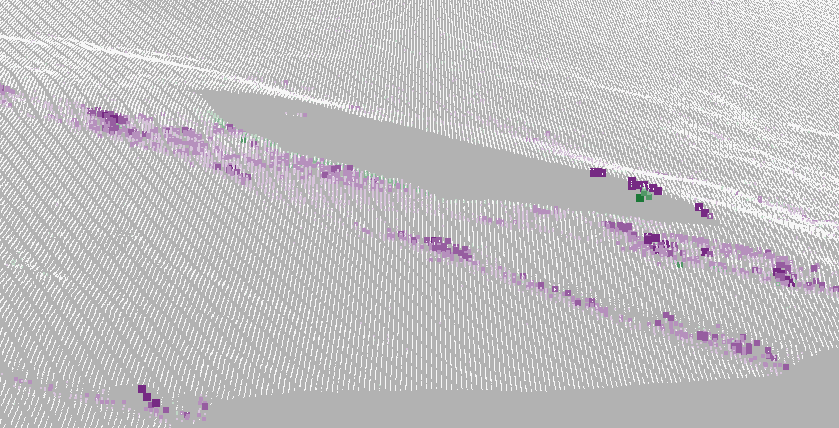} \\
&\vspace{0.1cm} \\
&\includegraphics[width=0.3\textwidth]{images2/raster05_20.png} \\
&\vspace{0.05cm} \\
(a) &(b)
\end{tabular}
\caption{(a) The point cloud coloured according to the distance to the 
raster with resolution 0.5 meters.
 (b) A detail.}
\label{fig:raster05_accuracy}
\end{figure}
The accuracy of a raster with resolution 0.5 meters is shown in Fig.~\ref{fig:raster05_accuracy}. The raster size is $1\,769 \times 1\,967$, the maximum distance between the point cloud and the 
raster surface is 2.96 meters and the average distance is 0.068 meters.
The accuracy is improved compared to the lower resolution raster, but the
distances between the surface and the point cloud are still larger than
for the LR B-spline surface. The largest distances are still to be found
in steep and ragged areas. The IDW interpolation method and the bi-linear
estimation of points in the surface are chosen for simplicity. Other 
interpolation approaches may result in better accuracy, but the effects of
the raster representation remain.

The raster surfaces are represented as GeoTIFF files with sizes 852KB and 14MB for the raster with resolution 2 meters and 0.5 meters, respectively.
The size of the LR B-spline surface is 2.5MB. GeoTIFF is a binary
format where the height values are stored as float, while the LR B-spline surface is stored in an ASCII file using doubles for 
the storage of coefficients. The height values of the raster with 2 meters resolution sum up to 217\,956 floats. The storage of the LR B-spline surface 
requires a total of 537\,175 numbers, of which 35\,801 are doubles and the 
remaining are integers. A compressed version of the LR B-spline
file has a size of 584KB.

Can  more iterations in the approximation algorithm improve the accuracy of the LR B-spline surface further? Increasing the number of iterations to
9, the number of points with a distance larger then 0.5 meters decreases
to 21\,446 while the number of coefficients grows to 203\,772 and the file size increases from 2.5MB to 16MB. The file size of the point cloud
is 359MB. After 16 iterations the least possible maximum distance of 1.19 meters is reached. The average distance is 0.052 meters and 476 points are further away from the surface than the given threshold. The surface size has grown 53MB. Continuing to 20 iterations another 5 points are within the threshold at the cost of increasing the surface size with 4MB. It is not necessarily a gain in insisting on the maximum possible accuracy. Accuracy must be balanced against the surface size and approximation of noise and outliers is not attractive.

\section{Export} \label{export}
The LR B-spline representation format is relatively new, and the current
support is limited. Thus, an option to export these surfaces to other formats
is crucial for the use of the LR B-spline format. Raster is the 
standard representation
for terrains and sea bed, and raster creation is a matter of regular 
evaluation, possibly with some adaption to features such as extremal points,
ridges and valleys. A raster representation will, in general, not support the 
same level of detail
as an LR B-spline representation of an area, but a high accuracy LR B-spline 
surface can give rise to rasters of
different resolution. Thus, the LR B-spline surface can serve as a master 
representation to be harvested according to needs. 

\begin{figure}
\centering
\begin{tabular}{ccc}
\includegraphics[width=0.3\textwidth]{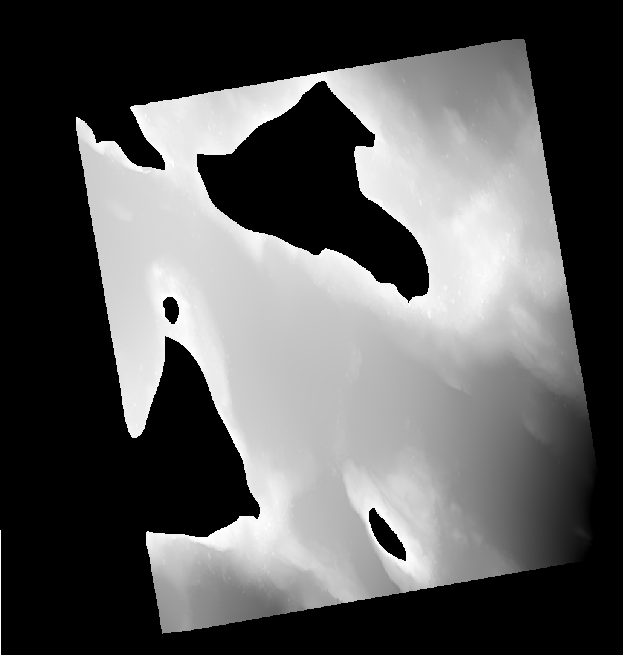}
&\includegraphics[width=0.3\textwidth]{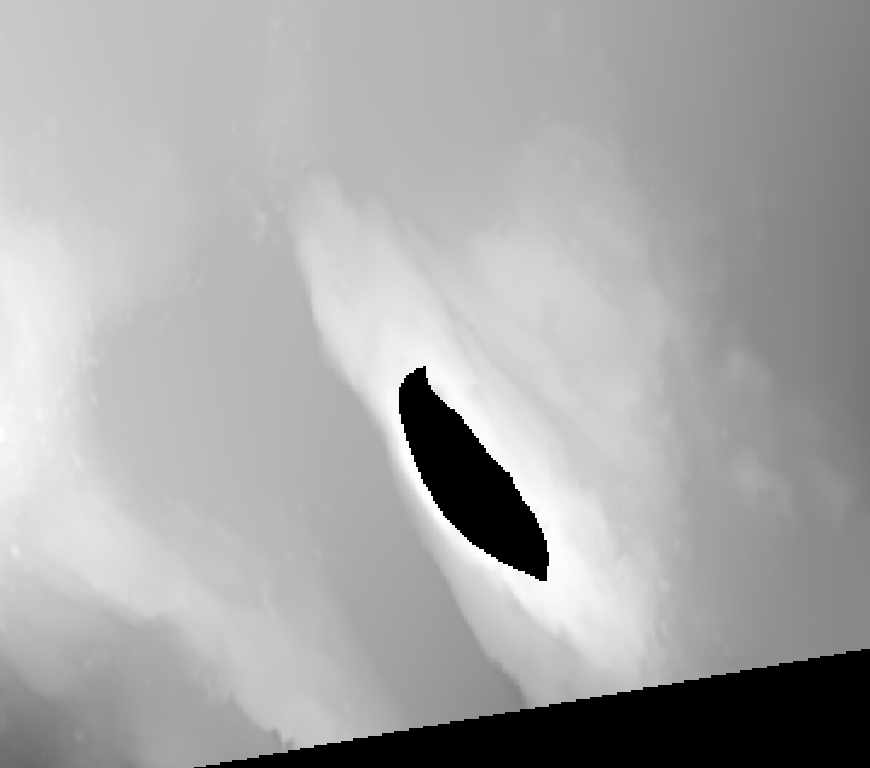}
&\includegraphics[width=0.3\textwidth]{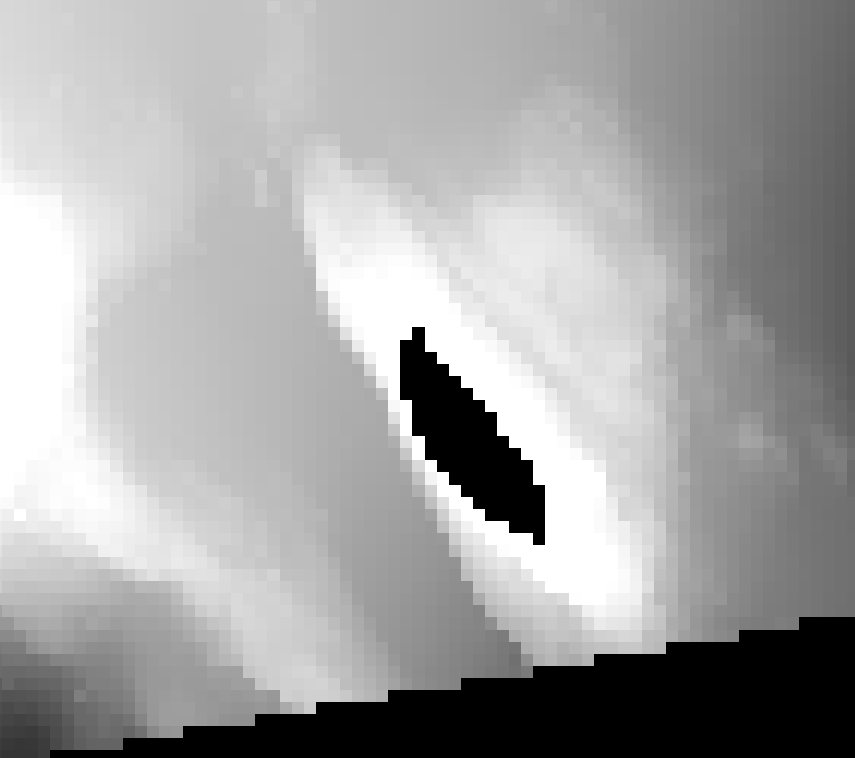}\\
(a) &(b) &(c)
\end{tabular}
\caption{(a) Raster with one meter resolution. (b) A detail. (c) Detail with
  five meters resolution}
\label{fig:DEM}
\end{figure}
Fig.~\ref{fig:DEM} shows raster representation at two different resolutions generated from
the LR B-spline surface in Fig.~\ref{fig:exs_approx}.

\subsection{Export as ISO 10303 LR B-spline surfaces}
ISO 10303 (STEP) is a standard for digital exchange of product data. It has
a broad scope, but CAD data exchange was an early adapter of the standard. 
Tensor product spline surfaces belong to the generic resources (Part 42) and in 
the context of the EC Factories of the Future project TERRIFIC (2011-2014), it was proposed 
to extend Part 42 with locally refined surfaces. The extension
was published as part of STEP in 2018. It is
ongoing work to make this format available for STEP users. In the future 
LR B-spline surfaces can be exchanged directly through STEP.

\subsection{Export as tensor product B-spline surfaces}
Direct export of LR B-spline surfaces through neutral exchange formats is an option
in the future.
Being a spline surface, an LR B-spline surface can be
expanded to a tensor product (TP) spline surface at the cost of a potentially
large increase in data size. This conversion contradicts the idea of locally 
refined splines. The size of the surface in 
Fig.~\ref{fig:exs_approx} would increase from 2.5MB to 17MB and the increase
in the number of polynomial elements is much higher due to different file 
formats. An LR B-spline surface could be exported as a set of Bezier surfaces, 
but as such a surface typically
contains a high number of polynomial elements, it is not obvious that this is a good 
solution. A better option is to represent the LR B-spline surface by a collection
of tensor product spline surfaces maintaining the feature of data size 
distributed according to needs.

\begin{figure}
\centering
\begin{tabular}{ccc}
\includegraphics[width=0.31\textwidth]{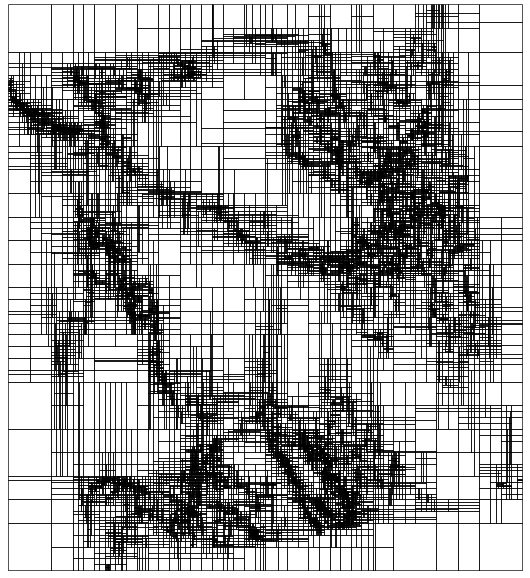}
&\includegraphics[width=0.31\textwidth]{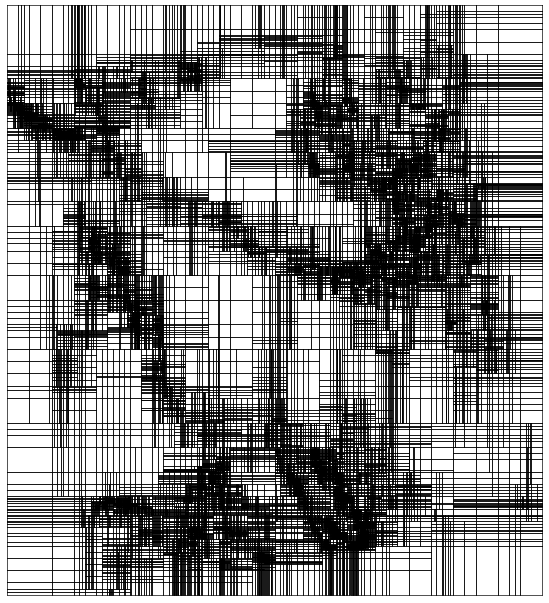}
&\includegraphics[width=0.31\textwidth]{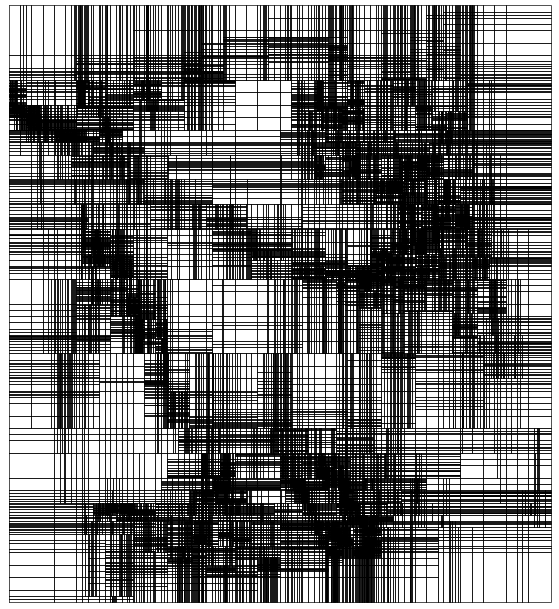}\\
(a) &(b) &(c) 
\end{tabular}
\caption{(a) LR B-spline surface mesh.  
(b) Mesh when the surface is divided into 666
TP surfaces. (c) Mesh when the surface is divided into 300
TP surfaces}
\label{fig:divideTP}
\end{figure}
Fig.~\ref{fig:divideTP} indicates how the LR B-spline surface shown in
Fig.~\ref{fig:exs_approx} can be divided into tensor product spline surfaces by
 the means of dedicated knot line insertions. Image (a) shows the polynomial patches of the LR B-spline surface,
(b) and (c) shows the patches after dividing the surface into 666 and 300
tensor product patches, respectively.

The division into TP surfaces is performed by a recursive algorithm. 
At each level it considers how
the current surface can be split by extending one knot line to cover the entire
surface domain. The candidate knot line must contain T-joints, that is at 
least one knot line in the other parameter direction must end at this knot line. The number of surface elements overlapping the knot line extension should be minimized and at the same time the knot line should 
divide the current surface into two surfaces with roughly the
same number of knots.
The balance between the two criteria varies throughout the recursion levels.
When an appropriate split is found, the algorithm proceeds to look for splits
in the two sub surfaces. The splitting algorithm stops when no sub surface
contains more segmented knot lines than a given threshold. Each sub
surface is expanded to a tensor product spline surface by adding missing
knot segments.

\section{Analysis Tools} \label{analysis}
\begin{figure}
\centering
\includegraphics[width=0.6\textwidth]{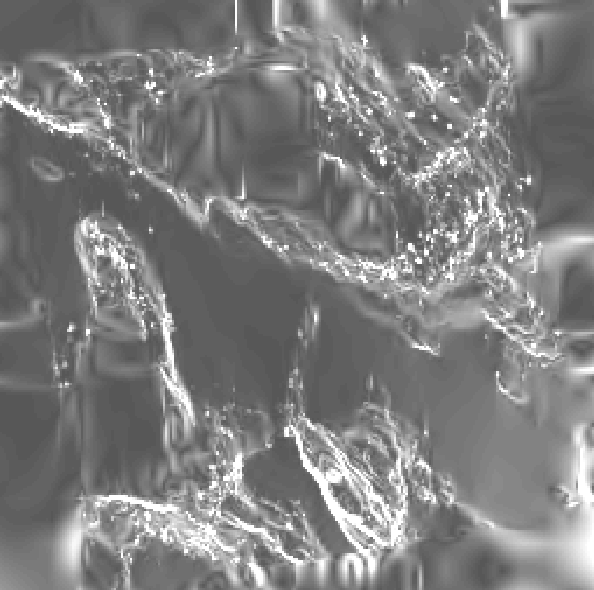}
\caption{The slope of the surface in Fig.~\ref{fig:exs_approx} (b) is evaluated in a raster}
\label{fig:slope}
\end{figure}
A GIS surface can be used to gain insight in the domain it represents. Along 
with visualization, properties such as slope, aspect and contour curves are
computed. Being a spline surface a point-wise computation of slope and
aspect is a straight forward task, an example of slope is shown in Fig.~\ref{fig:slope}. We will here go into some details regarding
contour lines and also present computation of minimum and maximum points.

\subsection{Contour lines} \label{contours}
Calculation of contour curves  is supported in GIS systems.
LR B-spline surfaces representing elevation are piecewise polynomial functions, and contour curves are curves
where the value of the spline function is constant. We want to find $f_a(t)=(f_1(t),f_2(t))^T \in R^2$ such that 
$F(f_1(t),f_2(t)) = a$ for an LR B-spline surface $F$ and an of elevation value $a$.

Instead of developing  interrogation functionality for LR B-spline surfaces, we split the LR B-spline surface into
a number of tensor product surfaces as described in the previous section. Then 
we  use the interrogation functionality of SINTEF's spline library, SISL~\cite{sisl} on each sub surface.

\begin{figure}
\centering
\begin{tabular}{cc}
\includegraphics[width=0.46\textwidth]{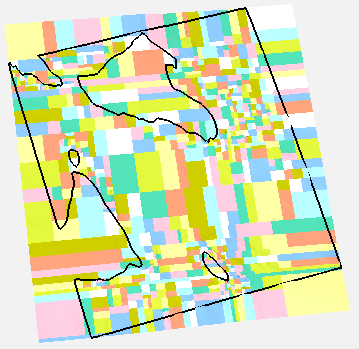}
&\includegraphics[width=0.49\textwidth]{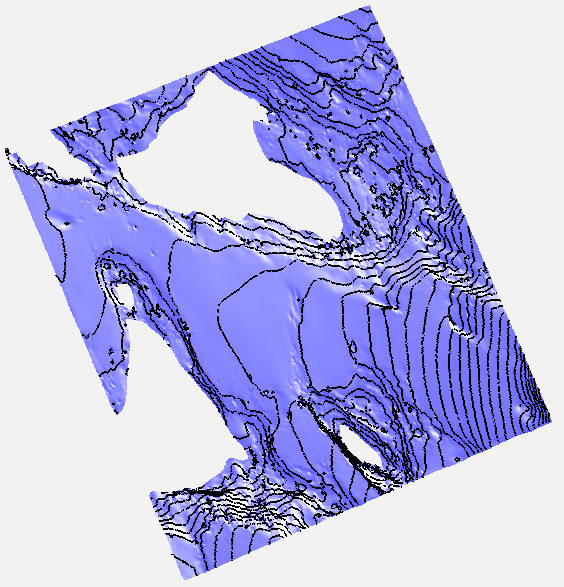}\\
(a) &(b)
\end{tabular}
\caption{(a) Division into TP surfaces. The trimming loops corresponding to the
LR B-spline surface are shown in black. (b) Surface with contours, one meter resolution}
\label{fig:contour1}
\end{figure}

Fig.~\ref{fig:contour1} (a)
shows how the example surface in Section~\ref{sec:example} is split into TP surfaces while (b) shows 
contour curves with one meter resolution. The contouring problem corresponds to computation of intersections between a parametric spline surface and an algebraic surface, a problem that is discussed in~\cite{intersect2}. The applied algorithm contains
several phases:
\begin{enumerate}
\item Divide the LR B-spline surface into a set of TP surfaces
\item For each value $a$ and each TP surface
\begin{enumerate}
\item Compute the topology of the contour curves using SISL. This is a recursive 
algorithm that finds a set of ``guide points'' on each curve branch.
\item Trace  each identified curve branch starting from an identified ``guide point''. Represent the curves traced out as  spline curves.
\end {enumerate}
\item For each value $a$: combine sub curves from different TP-surfaces into contour curves for the
entire LR B-spline surface.
\end{enumerate}

The contour curves will be accurate with respect to the given elevation values and
approximate the LR B-spline surface with respect to a given tolerance.

\begin{figure}
\centering
\begin{tabular}{cc}
\includegraphics[width=0.49\textwidth]{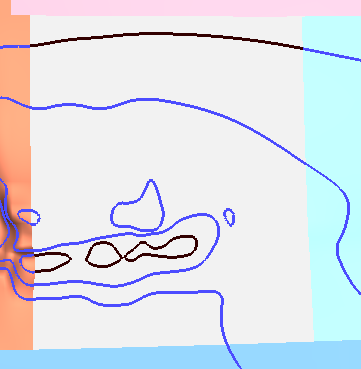}
&\includegraphics[width=0.45\textwidth]{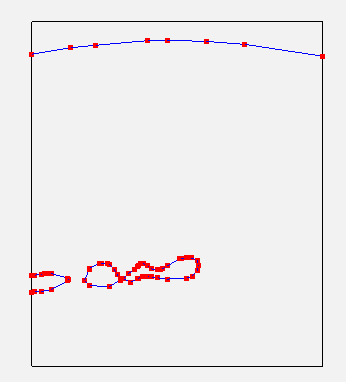} \\
(a) &(b) 
\end{tabular}
\caption{(a) Contour curves detail. 
(b) Result from the topology detection phase of computing the
black contour curve in (a).
Connected guide points describes 4 contour curves} \label{fig:contour2}
\end{figure}

An LR B-spline surface approximating an area with large shape variation will
contain many details, which again will lead to a complex pattern of contour
curves, see Fig.~\ref{fig:contour2} (a). The tensor product surfaces are distinguished by colour, and we will describe the computation
of a set of contour curves on the central TP surface
in some detail. 

The surface patch is in very shallow water and contains contour curves
at the depths of two, three and four meters. The pattern of contour curves at depth two meter (black curves in Fig.~\ref{fig:contour2}(a)) is complex with 
one curve passing between opposite boundaries, one curve with both endpoints at the left boundary and two very dense inner closed loops. The topology
detection is concerned with identifying this pattern and generating starting points for the proceeding tracing of the contour curves.
Also guide points that give a more complete description of the curves are computed. The algorithm uses a recursive approach: 
\begin{enumerate}
\item Check if there is any possibility that the elevation value $a$ is in the range of the
the current function. Otherwise, stop the computation.
\item Compute all intersection points at the surface boundaries. This is done recursively in the number of parameter
directions. The structure of the lower dimensional algorithm follows the pattern of the algorithm described here.
\item Check if there is any possibility that the elevation value is in the range of
the function excluding boundaries. If not stop the computation at
this recursion level. \label{recursion}
\item Check for a possible existence of inner closed intersection loops.
If not mark connections between intersection points at the boundaries and
stop the computation at the current recursion level.
\item Divide the surface into two or four.
\item Compute intersections between the given elevation value and the boundaries between sub surfaces.
\item For each sub surface go to~\ref{recursion}
\end{enumerate}
More details can be found in~\cite{recursive}.
Efficiency and robustness of the algorithm is reached through good 
interception methods and a clever strategy for dividing the surface into
sub surfaces. A discussion on subdivision strategies for surface intersections
can be found in~\cite{intersect}. A general rule is to subdivide at singularities and internal in closed loops.
A complex situation lead to more subdivisions and
consequently more guide points as shown in Fig.~\ref{fig:contour2} (b). The relatively simple curve at the top is identified with fewer subdivisions.

\begin{figure}
\centering
\includegraphics[width=0.5\textwidth]{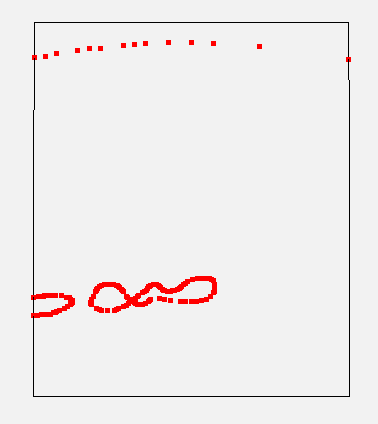} 
\caption{Tracing the result from guide points.}\label{fig:contour3}
\end{figure}
Consider again the current TP surface in Fig.~\ref{fig:contour2} and denote the surface $\tilde{F}(u,v)$. The surfaces is bi-quadratic and $C^1$ 
continuous. This degree of smoothness imply that the contour curves can turn
quite sharply at knot lines. The objective when tracing a contour curve is to 
describe the curve with sufficiently accuracy, handle sharp turns in the curve
and avoid jumping to a different contour curve. In the current case, the distance between the two closed loops is very small. The algorithm works as follows:
For a current point $(u_0,v_0)$ on the contour curve estimate the position of the next
point and iterate it onto a contour curve. The contour curve tangent direction at the point is, in the parameter domain, given by
$(-\frac{\partial F}{\partial v}(u_0,v_0), \frac{\partial F}{\partial u}(u_0,v_0))$  and the step length is deduced from the lengths of the first and second surface derivatives  at the point and the surface size. Ensure consistency of the directions
between the start point and the endpoint of the current segment, and between
both points and the midpoint of the segment. Other
intermediate points may also be involved in the verification of consistency. If  consistency is not verified, a new candidate
endpoint closer to the current one is computed and the process is repeated.
The tracing function uses information about the first and last guide point obtained from
the topology detector and whether or not the curve is closed. Typically, the
tracing function will produce a denser and more uniform sequence of points
than the topology detector, see Fig.~\ref{fig:contour3}. A simple configuration leads to less points than a more complex one.

\begin{figure}
\centering
\includegraphics[width=0.95\textwidth]{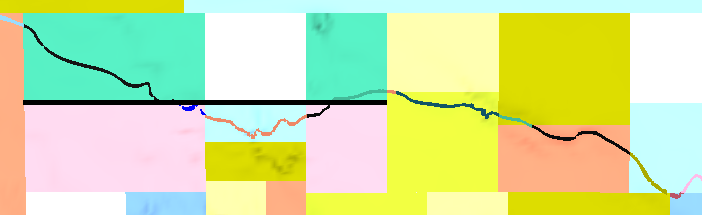} 
\caption{Joining curve fragments into complete contour curves across
  TP surface boundaries. Detail. The various curve segments have different colours and TP surfaces are also distinguished by colour.}
\label{fig:merge}
\end{figure}
The final step is, for each elevation value $a$, to join curve segments across the
boundaries of the tensor product surfaces, see Fig.~\ref{fig:merge}. Here
two pairs of curve segments belonging to an open curve need to be merged at the surface
interface marked with a horizontal black line. Note that there are still curve segments that need to be merged in the other parameter direction. During the processing of the 
interfaces between the TP
surfaces, more and more small curve fragments will be joined and
finally all contour curves will end at boundary curves of the
complete LR B-spline surface, they will be closed or meet at a branch point in the inner of the surface.

\subsection{Extremal points}
Information about minimum and maximum points (extremal points) can contribute 
to get a picture
of the behaviour of the terrain or sea bed in a given area. Such points are
also contained in the set of standard map information, provide
important information in the computation of a lower resolution representation of the
domain, and are candidates for landmarks. However, it is not
obvious how these points should be defined. The partial surface derivatives at an extremal 
point are zero and a maximum point has the least depth in the vicinity. Similarly,
is a minimum point situated at the locally greatest depth. This rules out saddle points
that typically will be placed between two extremal points of the same type.
Consider the surface visualized in Fig.~\ref{fig:exs_approx} (b) and Fig.~\ref{fig:contour1} (b).
This surface probably has one or a few global maxima, mostly placed 
at surface boundaries and also a few global minima. The surface
has a huge amount of local extrema and the higher the accuracy of the
approximation, the higher is the expected number of local extrema.

\begin{figure}
\centering
\begin{tabular}{cc}
\includegraphics[width=0.45\textwidth]{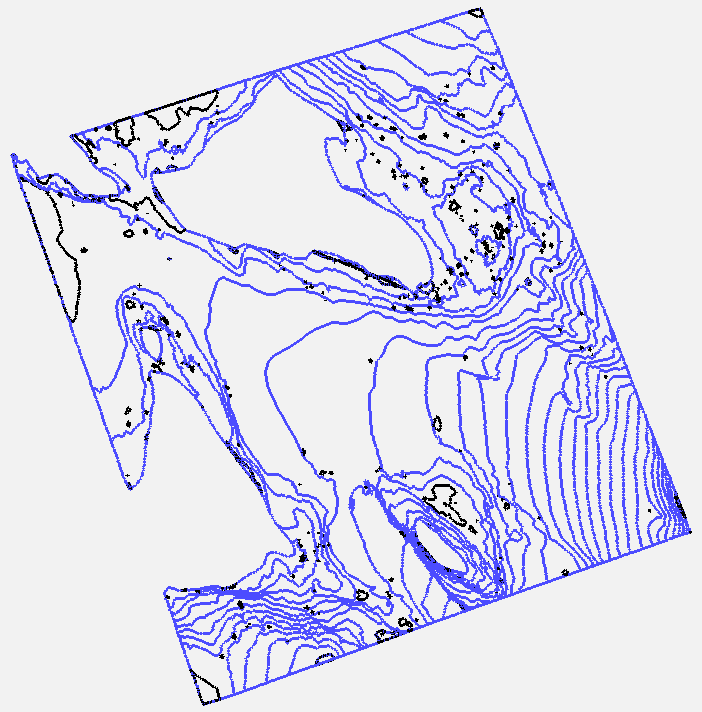}
&\includegraphics[width=0.45\textwidth]{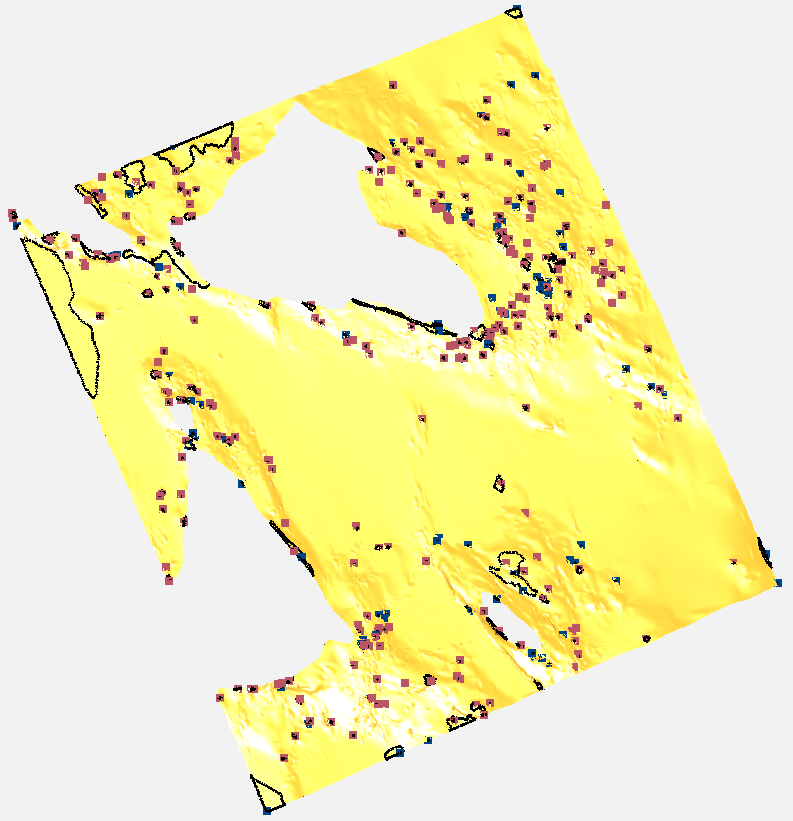} \\
\includegraphics[width=0.2\textwidth]{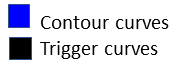}
&\includegraphics[width=0.25\textwidth]{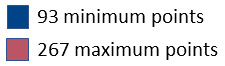} \\
(a) &(b) 
\end{tabular}
\caption{(a) Contour loops. The black curves trigger the computation of extremal points. 
  (b) The surface with trigger curves and identified extremal points. }
\label{fig:extpts}
\end{figure}

The granularity of the extremal point definition can be linked to the
granularity of the contour curves. Minimum and maximum points that do not 
deviate more from their surroundings than the interval between 
subsequent elevation values are not regarded as significant. Thus, a condition for the
existence of an extremal point in an area is that the extremal point is surrounded by
a closed contour curve or a contour curve meeting the surface boundaries.
Furthermore, there must be no contour curves inside an identified
one that are associated with the same type of extremal point.
If the terrain
height increases through a contour curve, before decreasing again to form a dump,
it can be appropriate to define both a minimum and a maximum point. 
Fig.~\ref{fig:extpts} (a) shows the contour curves of the test surface. The ones indicating an extremal point are black.  In
(b) the calculated extremal points are shown together with the surface and the
trigger curves. We see that the vast majority of extremal points are maximal 
and that they are quite densely distributed in some areas. Typically, the 
surface has local height variation in the vicinity of a contour level, which
can give extremal points of low significance. Thus,
post processing of the found extremal points will be needed to get a more
representative selection.

Each trigger curve defines an area in which one minimum or maximum
point is to be found. This area is represented as trimmed tensor product
surfaces using a similar approach as in the computation of contour curves.
Normally, one area leads to one TP surface. One global extremal point is
computed for each TP surface. The strategy is very similar to the topology
detection functionality for contour curves. Extremal points at the surface
boundaries are identified. A check for a possibly more extreme point in the
inner of the surface is performed. Since the surface is bounded by its coefficients   
the maximum value of the surface will be restricted by the maximum coefficient. If the control polygon given by the surface coefficients changes direction at most once in a parameter direction, then there will be at most one extremal point
in the surface, and an iteration
for the local extreme can be performed. Previously found, less extreme, points
are removed when a new extremal point is identified. If the current surface is too
complex to conclude with a result, the surface is subdivided and the search
continues on each sub surface.

The search for a global extremal point is performed on the local, non-trimmed 
tensor product surface. Thus, it is necessary to check whether the found
point is actually inside the trimmed surface. The terrain configuration combined with knowledge of an existing
extremal point indicates that this point often will lie in the trimmed surface. However, in areas with dense extremal points,
the algorithm may identify another extreme point lying within the non-trimmed TP surface. In our test case 52 out of 460 found extremal points
situated outside the trimmed surface. A fallback 
strategy is required to identify the extremal point inside the trimming loop.  A number of sample points
provide start points for an extremal point iteration. Then the most extreme point
is selected.

\FloatBarrier
\section{Adopting the approximation to emphasize particular areas} 
\label{adopt_apprx}
In shallow waters, as is the case in our example, an accuracy of the surface 
representation of in essence 0.5 meters may not be sufficient. In particular 
the boat traffic may need more exact information about shallows. We will, in this
section, investigate some strategies for ensuring sufficient clearance in critical
areas.
Table~\ref{tab:setup} summarizes the experiment setups investigated in this
section. The setup codes are constructed as follows: F means approximation with a fixed threshold and V with a variable threshold; S means that significant points are included in the computation; E indicates that a minimum element
size is defined; WM stands for the weighted mid surface; The first digit represents the maximum number of iterations and the
second, when given, the minimum element size.
Variable threshold and a restriction to the element size of the
surface is treated in Section~\ref{sec:var_threshold}, the weighted mid surface in
Section~\ref{sec:weighted} and significant points in Section~\ref{sec:significant}. 
\begin{table}
  \caption{Parameter settings for surface approximation. The threshold is given in
  meters, the number of iterations varies between seven and nine.}
\label{tab:setup} 
\footnotesize
  \begin{tabular} {llll}
\hline\noalign{\smallskip}
Setup & Threshold & No iter &  Comment\\
F7 & 0.5 & 7 & Example of Section~\ref{sec:example} \\
V7 & 0.20022 - 0.31176 & 7 & Variable threshold \\
V9 & 0.20022 - 0.31176 & 9 & Variable threshold \\
V9E1 & 0.20022 - 0.31176 & 9 & Restricting the element size to 1 $\times$ 1 square m. \\
V9E2 & 0.20022 - 0.31176 & 9 & Restricting the element size to 2 $\times$ 2 square m.\\
WM7 & 0.5 & 7 & Weighted mid surface \\
FS7 & 0.5 & 7 & 266 significant points, threshold 0.2 meters \\
FS9 & 0.5 & 9 & 266 significant points, threshold 0.2 meters \\
\noalign{\smallskip}\hline
\end{tabular}
\end{table}

\subsection{Surface approximation with a depth dependent threshold}\label{sec:var_threshold}
The specified accuracy threshold influences the approximation process in two ways: It governs where the 
 surface is refined; and, decides when the distance between
 the surface and the point cloud is sufficiently small. At each
 iteration step, the point cloud will be approximated as
 accurately as allowed by the spline space. A tighter
 threshold leads to more accurate surfaces, but it 
 also increases the surface size.
 
 The concept of local refinement allows more degrees of freedom
 in the surface where the point cloud has local details. A
 variable threshold can lead to a more accurate surface fitting 
 in critical areas and focus on representing the smooth
 component of the data in areas where the accuracy is less critical.
We now let the threshold depend on the sea depth at each point. Shallow water lead
to a stricter tolerance than great oceans depths.

\begin{table}
\caption{Figures of accuracy for approximations with varying number of
iterations and accuracy thresholds. The total number of points is 
11 150 110 and the size of the data file is 297 MB. 
For each parameter setup the sizes of the approximating surface in MB are reported
along with the number of surface coefficients and the maximum and average distances between the point set and
the surface. The number of points having a distance to the surface of less
than 0.2 meters, between 0.2 and 0.5 meters and more then 0.5 meters
are given in column 6, 7 and 8, respectively. The different
 setups are specified in Table~\ref{tab:setup}.}
\label{tab:accuracy} 
\footnotesize
\begin{tabular} {llllllll}
\hline\noalign{\smallskip}
Setup & Sf. size & No. coefs. & Max dist. & Av. dist. & $<$ 0.2
& $<$ 0.5 & $\ge$ 0.5 \\
\noalign{\smallskip}\hline\noalign{\smallskip}
F7 & 2.5MB & 33\,830 &2.85 m. & 0.068 m. & 10\,429\,550 & 685\,381 & 35\,179 \\
V7 & 16MB & 206\,295 &2.81 m. & 0.052 m. & 10\,753\,948  & 336\,828 & 27\,334 \\
V9 & 139MB & 1\,853\,313 &2.48 m. & 0.042 m. & 10\,957\,497 & 179\,908 & 12\,705 \\
V9E1 & 9.6MB & 139\,754 &2.81 m. & 0.055 m. & 10\,730\,787 & 388\,785 & 30\,538 \\
V9E2 & 3.7MB & 53\,824 &2.82 m. & 0.063 m. & 10\,496\,438 & 609\,920 & 43\,752 \\
\noalign{\smallskip}\hline
\end{tabular}
\end{table}
Table~\ref{tab:accuracy} compares the obtained surface approximation accuracy
for different refinement requirements. Setup F7 corresponds to the surface in
Fig.~\ref{fig:exs_accuracy}. 0.32\% of the points lie outside the given 
threshold of 0.5 meters. For the remaining cases having a varying threshold, 
the percentage of points with a distance larger than the threshold
is 2.68\%, 1.37\%, 3.12\% and 4.98\%, respectively. The data size
increases significantly with a stricter threshold and more iterations applied in
the approximation algorithm.
The element size of the final surface may in some areas be in the magnitude of the point cloud resolution.
To avoid approximation of noise or irrelevant details a restriction on the element size can be applied. This limits the surface size, but influences also
the approximation accuracy as can be seen in Table~\ref{tab:accuracy}, Setup V9E1 and V9E2.

\begin{figure}
\centering
\begin{tabular}{ccc}
\includegraphics[width=0.31\textwidth]{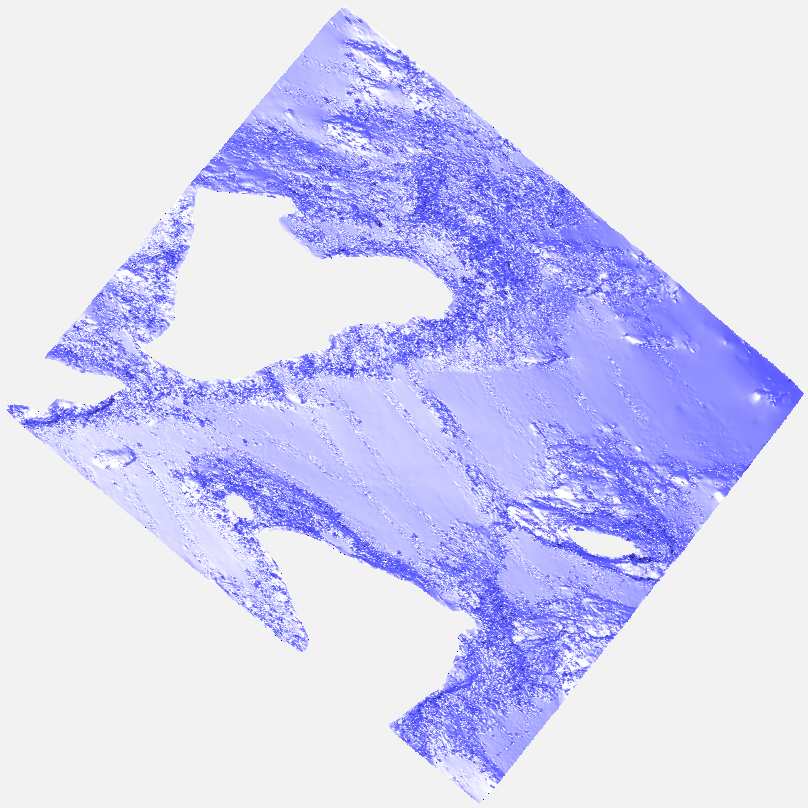}
&\includegraphics[width=0.29\textwidth]{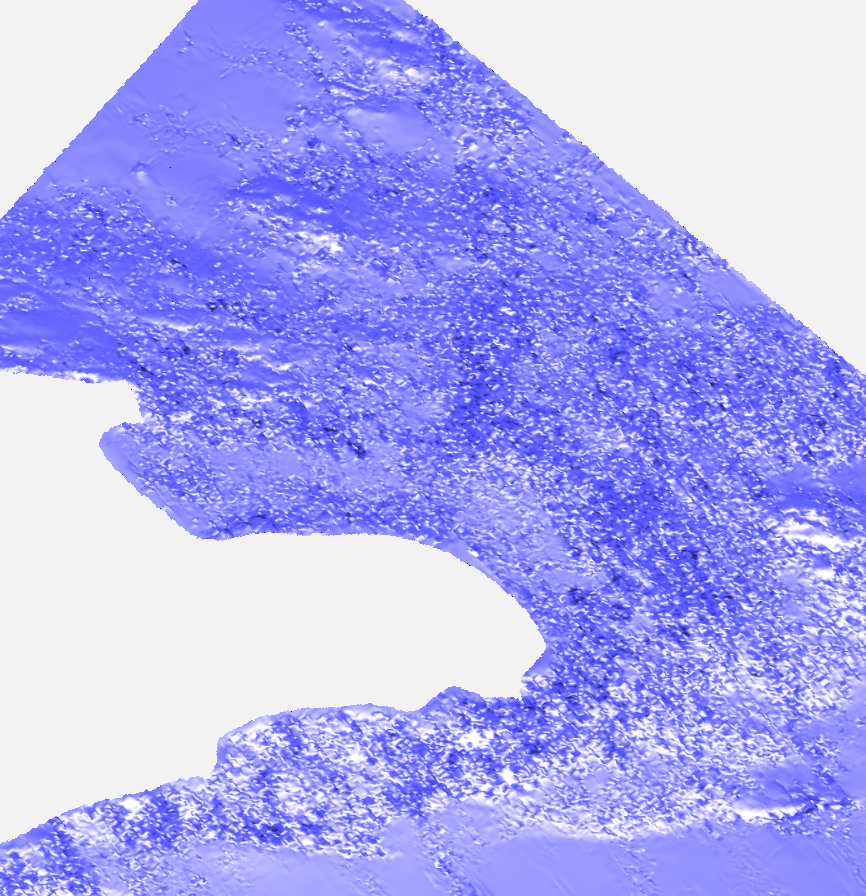} 
&\includegraphics[width=0.29\textwidth]{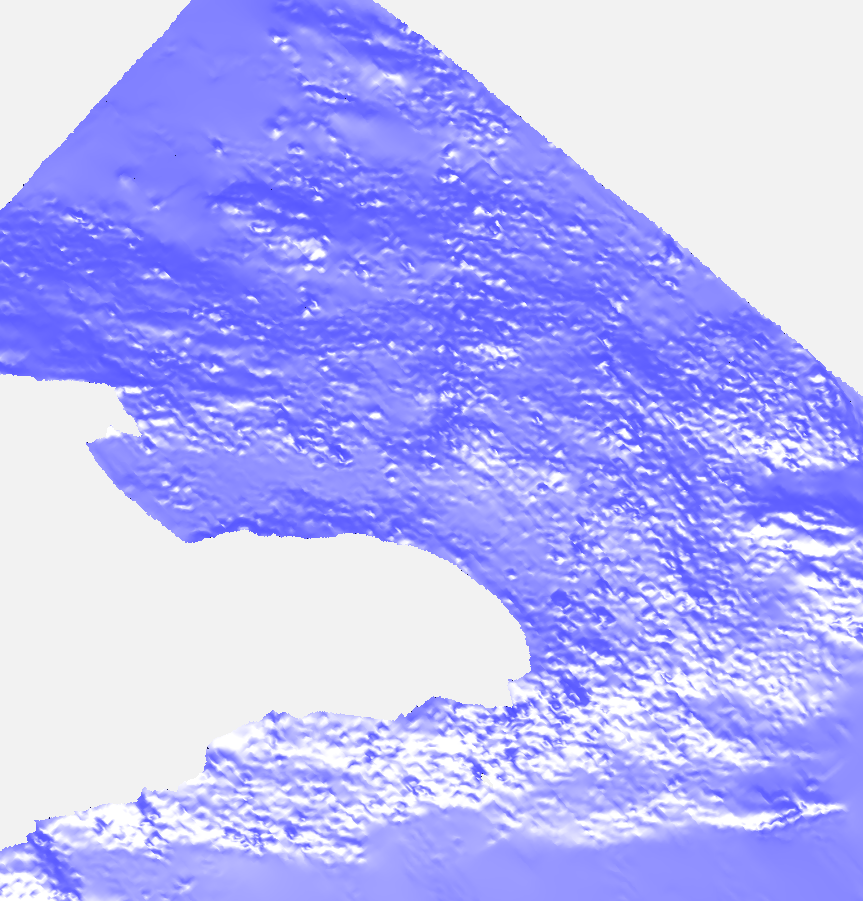} \\
(a) &(b) &(c)
\end{tabular}
\caption{(a) Approximating surface according to setup V9 (variable strict 
threshold, 9 iterations). (b) Detail of surface in (a). (c) Corresponding
detail for surface of setup V9E2 (as setup V9 with
element size restriction)}
\label{fig:approx02}
\end{figure}

The accuracy improves when the number of surface coefficients increases, but at a slower pace. The average distance between the surface and the points is
always kept low as the surface approximates smooth data very well. The
larger distances are caused by a non-smooth behaviour in the data that
can be due to
variation in the sea bed like stones, small inconsistencies in the data
set or even outliers, which is present in this data set. 
In fact, the surface tends to emphasize too much on
such variation if the degrees of freedom are high as can be seen in
Fig.~\ref{fig:approx02} where the surface with the highest number of coefficients (a) and (b) is very detailed in some areas. The less detailed surface in (c) also captures the main behaviour
of the sea bed, but omits some details. The approximation 
setup must reflect the data set properties and in this case the  point cloud contains noise.  Although being variable the provided threshold is strict and lead to a high number of
coefficients.
A depth dependent tolerance threshold is better suited to 
impose a less detailed approximation at large depths.

The threshold and the maximum number of iteration steps is important
for the approximation result. Currently these parameters are
set manually, but there is ongoing research to find a process related
stop criteria for the iteration and to set the threshold
with regard to point cloud properties.

The approximation threshold should not be so tight that we
focus on approximation of noise. However, there might be
areas where a failure to capture some details is severe. To create an overall smooth and
reasonably accurate surface with emphasis on particular details, there are
still opportunities with the LR B-spline approach.

\subsection{Limit surfaces} \label{sec:weighted}
The distances between a surface and the associated points vary depending on
the extent of steepness or
roughness in the sea bed/terrain and the point quality. In general, the distance will be larger
in non-smooth areas. The point cloud can be bounded by
limit surfaces, giving a means to enquire the uncertainty of the elevation 
represented by the given surface at every position in the
surface domain. The limit surfaces share the spline
space of the approximating surface and lies at opposite sides of the surface and the point cloud.

The initial point cloud is divided into two point clouds lying on either side 
of the approximating surface, from now on denoted the source surface. Each sub cloud give
rise to a residual set, one corresponding to points lying above the source surface and
one to points lying below. The residual sets are approximated individually by applying
several iterations of the MBA-algorithm. The resulting surfaces are defined with the
same collection of scaled tensor product B-splines as the source surface and no further refinement of the surfaces
is performed. The resulting residual surfaces approximate
the residual sets with high accuracy, but there might still
be some residuals that exceed the corresponding residual surface.

A post process ensures interpolation of the
remaining residuals. Given a residual $r_a$ at position 
$(x_a,y_a)$ where the distance $\Delta a = r_a - F_{upper}(x_a,y_a) > 0$. 
Here $f_{upper}$ is the residual surface corresponding to the
points lying above the source surface. Since the scaled tensor product B-splines
$N_B$ of $F_{upper}$ have the partition of unity property, we
know that
$$F_{upper}(x_a,y_a)+\Delta a = \sum_{B\in \mathcal{B}_j} (c_B + \Delta a) N_B(x_a,y_a)$$ following the notation of Equation~\ref{eq:LRSurf}.
$\mathcal{B}_j$ is the collection of B-splines corresponding to of the residual surface.
For each $B\in \mathcal{B}_j$ we add the maximum distance between the
residuals in the B-spline support and the surface to the
corresponding coefficients. This ensures that $F_{upper}(x_a,y_a) \ge r_a$ for all residuals. The lower limit surface is treated
similarly.
Finally, the two limit surfaces are created by adding 
the residual surfaces to the source surface.

\begin{figure}
\centering
\begin{tabular}{cc}
\includegraphics[width=0.45\textwidth]{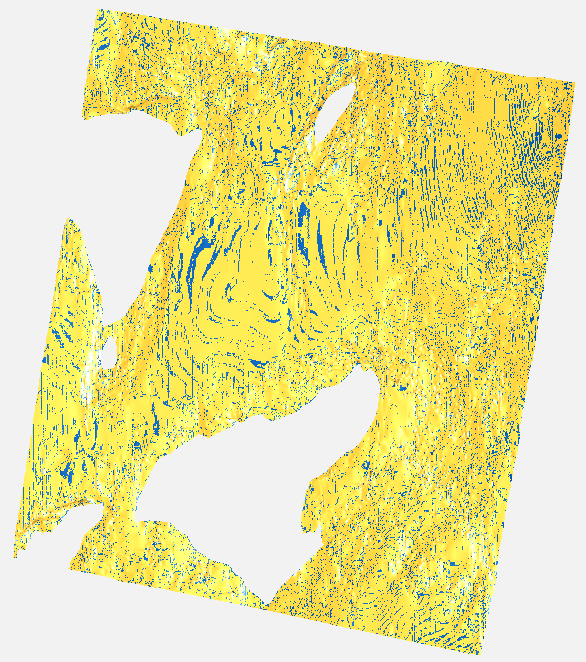}
&\includegraphics[width=0.45\textwidth]{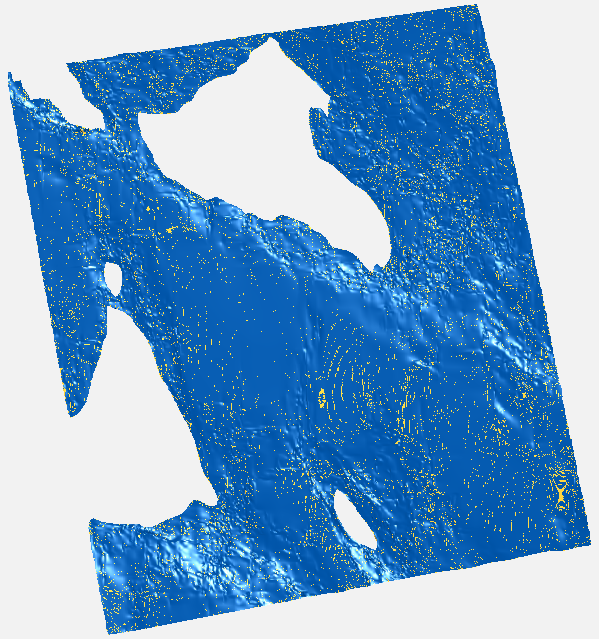} \\
\includegraphics[width=0.2\textwidth]{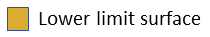}
&\includegraphics[width=0.2\textwidth]{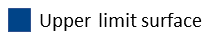} \\
(a) &(b) 
\end{tabular}
\caption{(a) Limit surfaces seen from below.
(b) Limit surfaces seen from above.} 
\label{fig:limit}
\end{figure}

An example of limit surfaces is shown in
Fig.~\ref{fig:limit}. The source surface is
generated by Setup F7 in Tables~\ref{tab:setup} and~\ref{tab:accuracy}. The maximum distance
between the two limit surfaces is 4.26 meters, the minimum distance is zero,
and the average distance is 0.93 meters. Seen from below, Fig.~\ref{fig:limit}(a), the
lower limit surface is dominant, but the upper limit surface can
be glimpsed where the surfaces are tight. Similarly is the upper
limit surface dominant when seen from above (b).
\begin{figure}
\centering
\begin{tabular}{cc}
\includegraphics[width=0.45\textwidth]{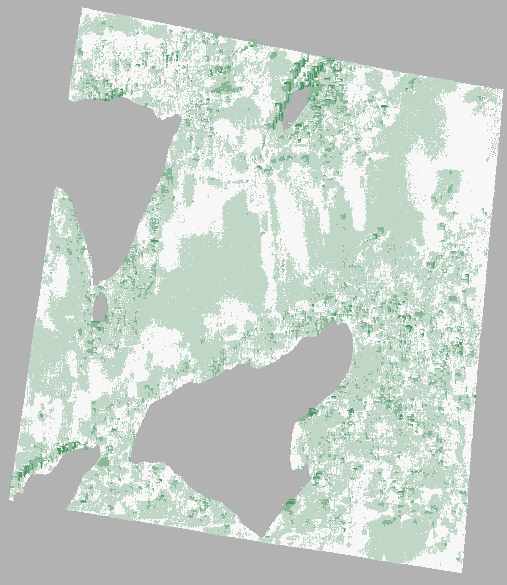}
&\includegraphics[width=0.47\textwidth]{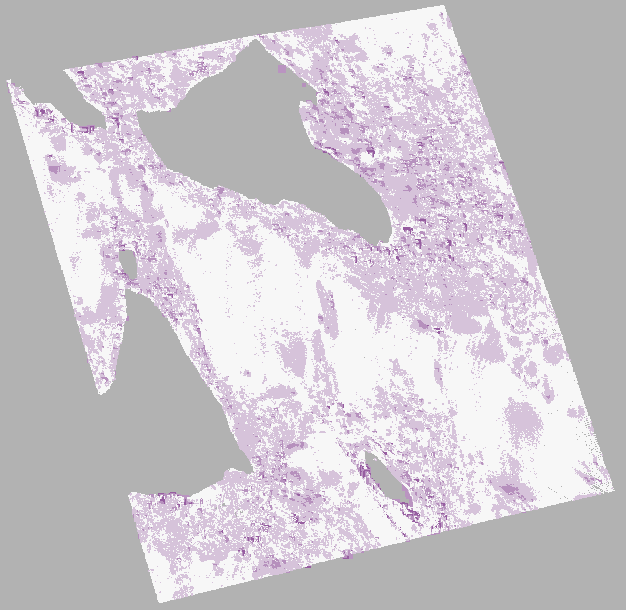} \\
\includegraphics[width=0.32\textwidth]{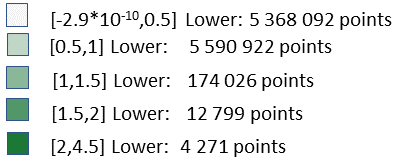}
&\includegraphics[width=0.32\textwidth]{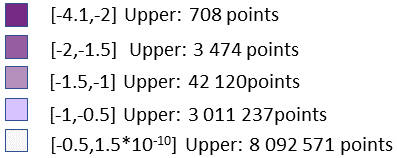}\\
(a) &(b) 
\end{tabular}
\caption{
(a) Point cloud coloured according to the distance to the lower limit surface. 
(b) Point cloud coloured according to the distance to the upper limit surface. }
\label{fig:limit2}
\end{figure}
In Fig.~\ref{fig:limit2} the distances between the
original point cloud and the limit surfaces are visualized (a) and (b). The points are coloured according to the distances to the
associated surface and the colour scale and how the points are
distributed with regard to distance are shown together with
the surfaces. The distance
range between the point cloud and the upper limit surface is
$[-4.06m.,1.5\times10^{-10}m.]$ and the average distance is 0.41
meters. The corresponding numbers for the lower limit surface
are $[-2.0\times10^{-10}m.,4.48m.]$ and 0.53 meters.

The upper limit surface is guaranteed to lie above all measurement points.
Thus, it can serve as a conservative model with respect to the sea depth
at shallows. 
A weighted surface between the source surface and the upper limit
surface can be created to ensure a safety zone at shallows
while restoring to the source surface at greater depths. Let $F_w$ be the weighted surface, $F_{source}$ the
source surface and $F_{upper}$ the upper limit surface. The surfaces correspond to
the same collection of B-splines. Following Equation~\ref{eq:LRSurf} we have
\begin{eqnarray*}
F_w(u,v) = &\alpha(u,v)F_{source}(u,v)+(1-\alpha(u_B,v_B))F_{upper}(u,v) = \\
&\sum_{B\in \mathcal{B}_j} \big(\alpha(u_B,v_B) c_{source,B}+(1-\alpha(u_B,v_B)) c_{upper,B}\big)N_{B}(u,v)
\end{eqnarray*}
where the blending factor $\alpha$ depends on the elevation in a
given point $(u,v)$ and a given elevation range $[d_1,d_2]$ where the two
surfaces are blended. The scaled B-spline $N_B$ is a scaled tensor product B-splines and the parameter value, $(u_B,v_B)$, of
the coefficient corresponding to the B-spline is
given by the knot vectors $\mathbf{u}=(u_0,\ldots,u_{p_1+1})$ and $\mathbf{v}=(v_0,\ldots,v_{p_2+1})$ of the B-spline as the Greville point
$$(u_B,v_B) = \Big( \sum_{j=1}^{p_1}\frac{u_j}{p_1},\sum_{k=1}^{p_2}\frac{v_k}{p_2}\Big).$$ Let $d=0.5*(F_{source}(u_B,v_B)+F_{upper}(u_B,v_B))$. Now, we can define $\alpha$ as
$$
\alpha(u_B,v_B) = \begin{cases} 1; &\text{if } d \le d_1, \\
0; & \text{if } d \ge d_2, \\
(d_2-d)/(d_2-d_1); &\text{otherwise}. 
\end{cases}
$$

Given $d_2=0$ and $d_1$ at a specified depth  the weighted mid-surface will be close to the
upper limit surface in shallow areas and coincide with the original approximating
surface in deeper water. 
\begin{figure}
\centering
\begin{tabular}{ccc}
\includegraphics[width=0.33\textwidth]{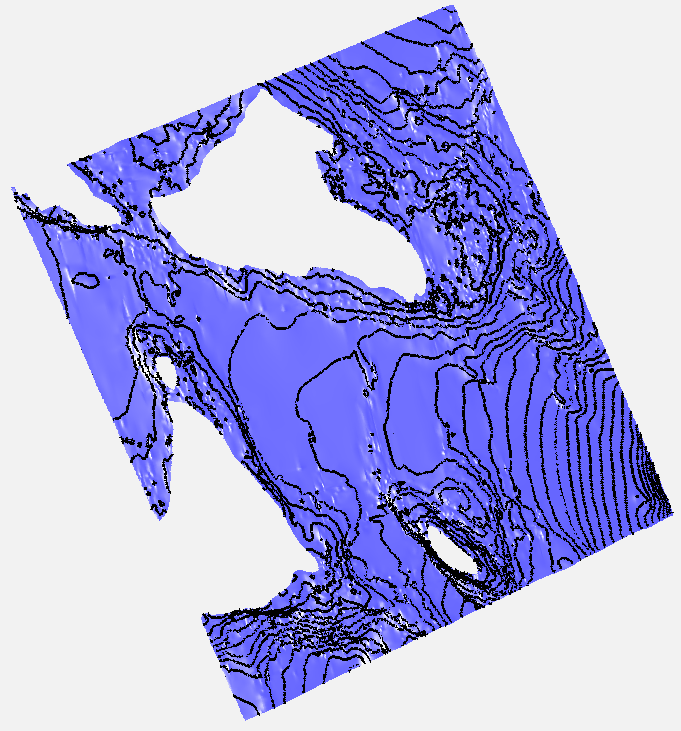}
&\includegraphics[width=0.33\textwidth]{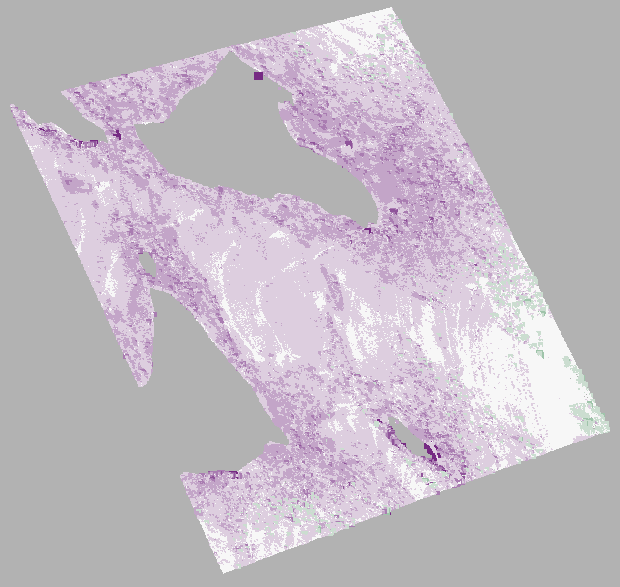}
&\includegraphics[width=0.26\textwidth]{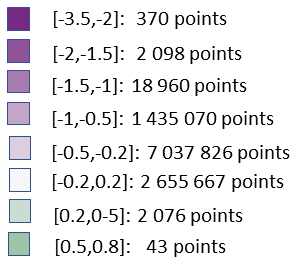} \\
(a) &(b) 
\end{tabular}
\caption{(a) Weighted mid surface with contour curves. 
(b) Data points coloured according to distance to weighted mid surface. The colour scale and number of points in each colour range is given.}
\label{fig:weighted}
\end{figure}
A weighted middle surface for the approximating surface of Setup F7 and the 
corresponding upper limit surface is shown in Fig.~\ref{fig:weighted}(a) along with contour curves with  resolution 1 meter. Here the limits of the transition zone are $d_1=-20$ and $d_2=0$. The contour curves differ slightly from the curves in Fig.~\ref{fig:contour1}, but the
main pattern is the same. As seen in
(b), the majority of the points lie below the surface. 2\,119 points lie above
the surface with a distance larger than 0.2 meters and these points are
situated at the largest depths ~in the data set. The maximum distances
between the surface and the points are 0.76 meters above the surface and
3.22 meters below. The average distance is 0.32 meters. 

\subsection{Approximation with significant points} \label{sec:significant}
The weighted mid-surface construction can be vulnerable
for outliers as the limit surface interpolates them. An approach for enforcing sufficient accuracy at shallows while maintaining a good overall accuracy is the concept of significant points.

Both least squares approximation and the MBA algorithm are able to 
weigh points individually. The surface will approximate points with
a higher weight more closely and it is possible to do additional refinement in areas with significant points. In this construction the weight associated with significant points is five times the
weight of ordinary points. Furthermore, an extra final iteration with the MBA method where significant points are weighted 50 times the weight of other points is applied if the
accuracy at significant points does not satisfy the threshold. Thus, defining significant
points at critical areas in shallow water and setting a stricter threshold for these points can address the need for a very
accurate depth representation at shallows by a lean and smooth surface.

\begin{figure}
\centering
\begin{tabular}{cc}
\includegraphics[width=0.45\textwidth]{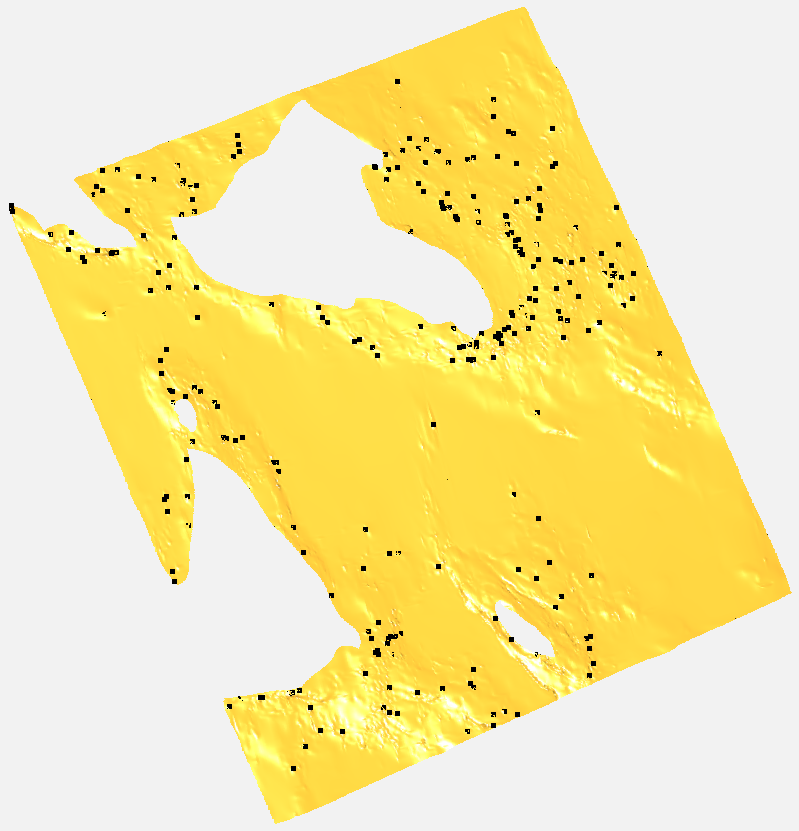}
&\includegraphics[width=0.45\textwidth]{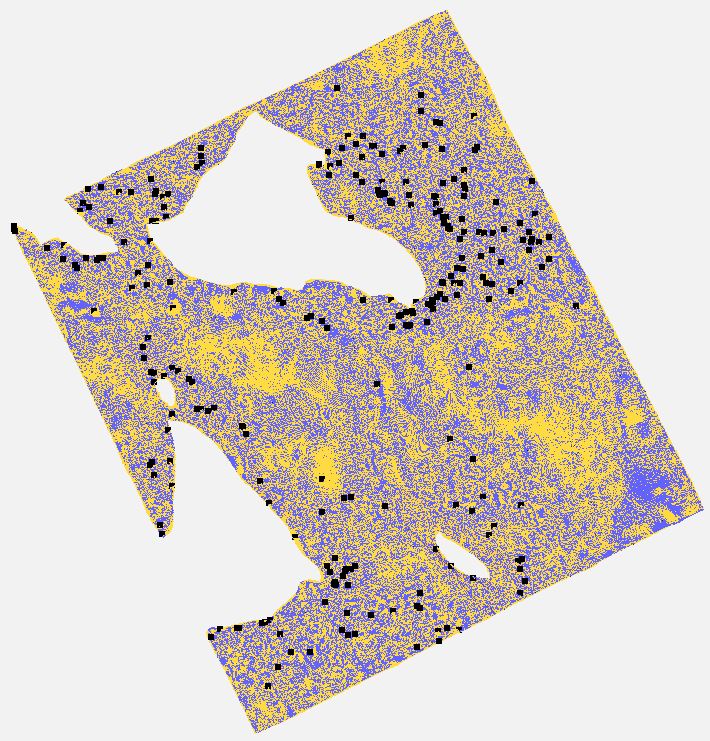} \\
\includegraphics[width=0.405\textwidth]{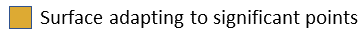}
&\includegraphics[width=0.46\textwidth]{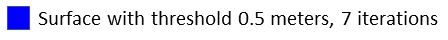}
\\
(a) &(b) 
\end{tabular}
\caption{(a) Significant points with surface approximation. The threshold corresponding to the complete data set is 0.5 meters, while the
threshold for the significant points is 0.2 meters (Setup FS7). (b) Surfaces from setup F7 and FS7. The significant points are included.}
\label{fig:signpts1}
\end{figure}

The election of significant points is crucial to get the wanted result:
good accuracy in critical points, high safety at shallows, low data size 
and avoiding  modelling of noise. We have selected 266 points from the original
data set. The points lie in the vicinity of the extremal points shown
in Fig.~\ref{fig:extpts} (b). One data point close to each  maximum point
having a depth of less than ten meters, is selected as a significant point. The resulting points have the minimum depths of all data points in their neighbourhood. Fig.~\ref{fig:signpts1} (a) shows the resulting surface along with the significant points. In (b) also the surface from Setup F7 is included. The two surfaces are roughly
similar, but differ slightly in some places. It is a tendency that
the surface including approximation of significant points lie above the other
in areas where the sea bed is smooth. 

\begin{figure}
\centering
\footnotesize
\begin{tabular}{cc}
\includegraphics[width=0.45\textwidth]{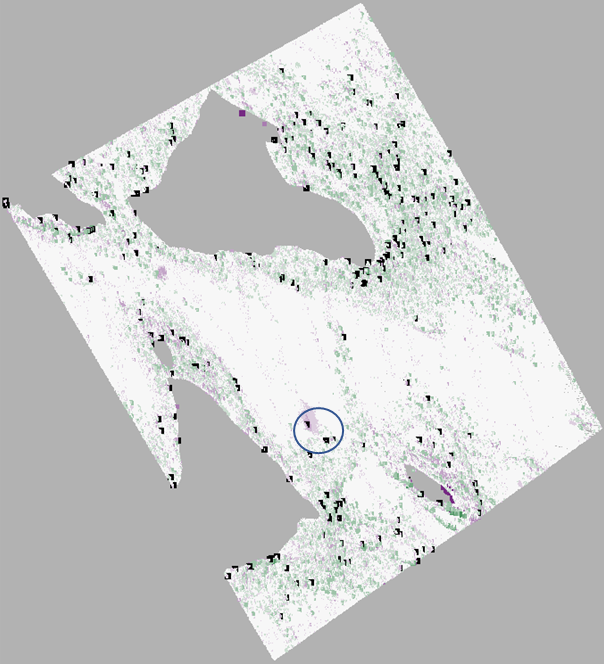}
&\includegraphics[width=0.47\textwidth]{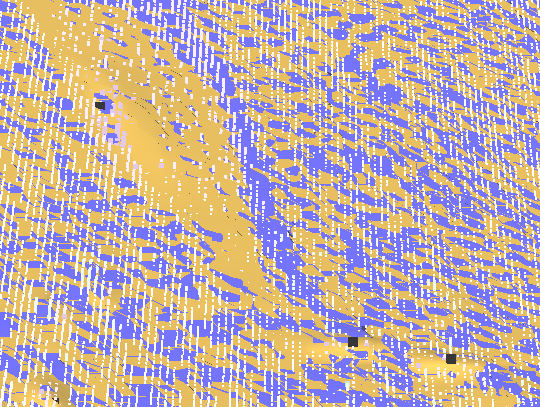}\\
\includegraphics[width=0.3\textwidth]{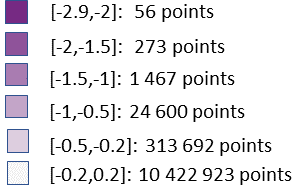}
&\includegraphics[width=0.3\textwidth]{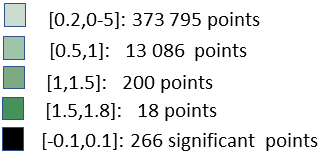}\\
(a) & (b) 
\end{tabular}
\caption{(a) Initial point cloud coloured according to the distance to the approximating surface where significant points have a stricter threshold. (b) Detail with both surfaces, the significant
points and the initial point cloud in the area marked with s
circle in (a).}
\label{fig:signpts2}
\end{figure}
Fig.~\ref{fig:signpts2} (a) shows the original point cloud coloured with respect to the
distance to the surface fitting also significant points (Setup FS7) and the significant
points. An area is marked for further study. The corresponding part of the two
surfaces, the point cloud coloured according to distance and
the significant points are shown in (b). The significant
points tend to reduce the depth represented by the corresponding
surface in an area in the vicinity of such a point.

\begin{table}
\caption{Figures of accuracy for approximations with different setups regarding
the number of iterations, tolerance threshold, significant points and
weighted mid surface construction. The three first cases are included from 
Table~\ref{tab:accuracy} as references.}
\label{tab:accuracy2}
\begin{tabular} {llllllll}
\hline\noalign{\smallskip}
Setup & Sf. size & No. coefs. & Max dist. & Av. dist. & $<$ 0.2
& $<$ 0.5 & $\ge$ 0.5 \\
\noalign{\smallskip}\hline\noalign{\smallskip}
F7 & 2.5MB & 33\,830 & 2.85 m. & 0.068 m. & 10\,429\,550 & 685\,381 & 35\,179 \\
V7 & 16MB & 206\,295 & 2.81 m. & 0.052 m. & 10\,753\,948  & 336\,828 & 27\,334 \\
V9 & 139MB & 1\,830\,313 & 2.48 m. & 0.042 m. & 10\,957\,497 & 179\,908 & 12\,705 \\
WM7 & 2.5MB & 33\,830 & 3.49 m. & 0.32 m. & 2\,655\,667 & 7\,037\,902 & 1\,456\,541 \\
FS7 & 2.6MB & 35\,491 & 2.9 m. & 0.068 m. & 10\,422\,923 & 687\,487 & 39\,700 \\
FS9 & 17MB & 217\,934 & 2.47 m. & 0.057 m. & 10\,669\,437 & 459\,456 & 21\,217 \\
\noalign{\smallskip}\hline
\end{tabular}
\end{table}

Table~\ref{tab:accuracy2} relates the weighted middle surface construction and approximation
with significant points to pure point approximation. The weighted mid surface has lower accuracy than the other
approaches. This is expected. The surface is lifted upwards to obtain a
safety depth at shallows. The accuracy when significant points 
are included compares quite equitable to approximations without significant
points, but the ordinary points are slightly punished for the extra emphasize on the significant points. For Setup FS7 and FS9 significant points with an associated tolerance of
0.2 are given. 
The maximum and average distance to the significant points after 7 iterations 
were 0.099 meters and 0.0008 meters, respectively, and the maximum distance between the surface and
the points lying above the surface is 1.8 meters. After 9 iterations the maximum distance between the
surface and the significant points is 0.023 meters and the average distance is 0.002 meters. The maximum distance
to points above the surface is reduced to 1.68 meters.
After 7 iterations, a final surface approximation with increased weight on the significant
points was applied. This was not required using 9 iterations iteration steps.

\section{Conclusion and further work} \label{conclusion}
The LR B-spline surface format is a new representation format that provides
a middle road between the rigid, but effective regularity of the raster format
and the large flexibility of triangulated surfaces. LR B-spline surfaces are
smooth and can, due to adaptivity, represent local detail without a drastic increase in data
size. Being a novel format, it is not supported in GIS systems, but
the surfaces can be exported as rasters in various resolutions as well as
collections of tensor product spline surfaces. 

Having a spline format, LR B-spline surfaces are well suited for
computation of various properties. We have
discussed the computation of contour curves and extremal points. The latter
will require some post processing to get a clean collection of points. This
will be looked into in the future. We will also look at a representation of
slope and aspect that gives further insight in the information these
quantities possess. We are for instance interested in a kind of contour curves
for slope. 

Generalization is the process of computing a lower resolution terrain or
seabed model from a high resolution one. Advanced methods take features
such as extremal points, ridges and valleys into account. Generalization is
not yet investigated for LR B-spline representations, but it is an
important and interesting topic for further work.

Outliers may obstruct the surface computation. Single points do not influence
the surface significantly, but collection of outliers points and outliers in 
areas with sparse data points will drag the surface in their direction. 
Moreover, the existence of outliers will obstruct the accuracy information.
To some extent outliers can be identified during the computation of the
approximating surface and removed during the iterative algorithm presented in
Section~\ref{sec:LR}. However, this approach is mostly suited for single
outlier points and there is always a risk of removing significant points. The
outlier problem is an important topic for further work.

Ongoing research focuses on the approximation threshold and a
stop criterion for the number of iterations in the approximation
algorithm. The intent is to find a balance between the accuracy of the approximation and
the number of coefficients representing the surface. The idea
is to get an accurate surface that do not adapt to noise in the
data.

 \vskip 0.5cm
{\noindent
 {\bf Data availability}
The data set used in the computations belong to an area in Norway called S\o re Sunnm\o re
where data is made generally available. The actual data set is from an older acquisition and
are available on request.
}

\end{document}